\documentclass[draft]{article} 
\usepackage{amssymb,amsmath,amscd,amsthm} 

\usepackage[all]{xy} 
\usepackage{color}

\font\tenrsf=rsfs10 at 11pt
\font\sevenrsf=rsfs7 at 8pt
\font\fiversf=rsfs5 at 6pt
\newfam\rsffam
\textfont\rsffam=\tenrsf
\scriptfont\rsffam=\sevenrsf
\scriptscriptfont\rsffam=\fiversf
\def\rond#1{{\tenrsf\fam\rsffam#1}}

\renewcommand{\theequation}{\arabic{section}.\arabic{equation}}
\newtheorem{theorem}{Theorem}[section]
\newtheorem{lemma}[theorem]{Lemma}
\newtheorem{proposition}[theorem]{Proposition}
\newtheorem{corollary}[theorem]{Corollary}

\theoremstyle{definition} 
\newtheorem{definition}[theorem]{Definition}
\newtheorem{example}[theorem]{Example}
\theoremstyle{remark} 
\newtheorem{remark}[theorem]{Remark} 

\newcommand{\cA}{\mathcal{A}}
\newcommand{\cB}{\mathcal{B}}
\newcommand{\cC}{\mathcal{C}}
\newcommand{\cI}{\mathcal{I}}
\newcommand{\cJ}{\mathcal{J}}
\newcommand{\cG}{\mathcal{G}}
\newcommand{\cH}{{}^c\!H}
\newcommand{\cL}{\mathcal{L}}

\newcommand{\cN}{\mathcal{N}}

\newcommand{\cV}{\mathcal{V}}
\newcommand{\dA}{d_A}
\newcommand{\Dr}{\rond{D}}
\newcommand{\Fr}{\rond{F}}
\newcommand{\Gr}{\rond{G}}
\newcommand{\Hr}{\rond{H}}

\newcommand{\Kr}{\rond{K}}
\newcommand{\Lr}{\rond{L}}

\newcommand{\Qr}{\rond{Q}}
\newcommand{\Sr}{\rond{S}}
\newcommand{\Bc}{\mathcal{B}}
\newcommand{\Dc}{\mathcal{D}}
\newcommand{\Kc}{\mathcal{K}}

\newcommand{\ccV}{{}^c\cV}
\newcommand{\comp}{\mathrm{comp}}
\newcommand{\ctX}{{}^cT\oX}
\newcommand{\ctsX}{{}^cT^*\oX}
\newcommand{\cun}{\cC^{\infty}}
\newcommand{\cunc}{\cC^{\infty}_c}
\newcommand{\cz}{\mathbb{C}}
\newcommand{\Diff}{\mathrm{Diff}}
\newcommand{\Diffc}{\mathrm{Diff}_c}

\newcommand{\dR}{{\rm dR}}
\newcommand{\dtt}{{d_{\tit}}}

\newcommand{\ess}{\mathrm{ess}}
\newcommand{\Ext}{\mathrm{Ext}}
\newcommand{\Hom}{\mathrm{Hom}}

\newcommand{\Krke}{\Kr_{k,\varepsilon}}
\newcommand{\loc}{\mathrm{loc}}
\newcommand{\oX}{\overline{X}}
\newcommand{\px}{\partial_x}
\newcommand{\rz}{\mathbb{R}}
\newcommand{\supp}{\mathrm{supp}}

\newcommand{\R}{\mathbb{R}}
\newcommand{\N}{\mathbb{N}}
\newcommand{\qz}{\mathbb{Q}}
\newcommand{\Z}{\mathbb{Z}}
\newcommand{\siges}{\sigma_\ess}

\newcommand{\tit}{\tilde{\theta}}

\newcommand{\tM}{\tilde{M}}

\newcommand{\ucz}{\underline{\cz}}
\newcommand{\vol}{\operatorname{vol}}
\newcommand{\Vol}{\operatorname{Vol}}

\newcommand{\zz}{\mathbb{Z}}
\newcommand{\slim}{\mathrm{s-}\!\lim} 
\def\cchi{\raisebox{.45 ex}{$\chi$}}
\def\build#1_#2^#3{\mathrel{\mathop{\kern 0pt#1}\limits_{#2}^{#3}}}
\def\pp{<}
\def\pg{>}
\begin{document}
\title{Spectral analysis of magnetic Laplacians on conformally
cusp manifolds}
\author{Sylvain Gol\'enia and Sergiu Moroianu}
\date{\today}
\maketitle

\begin{abstract}
We consider an open manifold which is the interior of a compact manifold
with boundary. Assuming gauge invariance, we classify magnetic fields 
with compact support into being trapping or non-trapping. 
We study spectral properties of the associated magnetic Laplacian for a 
class of Riemannian metrics which includes complete hyperbolic metrics of 
finite volume.
When $B$ is non-trapping, the magnetic Laplacian has nonempty
essential spectrum. Using Mourre theory, we show the absence of
singular continuous spectrum and the local finiteness of the point
spectrum. When $B$ is trapping, the spectrum is discrete and obeys
the Weyl law. The existence of trapping magnetic fields with compact
support depends on cohomological conditions, indicating a new and
very strong long-range effect. 

In the non-gauge invariant case, we exhibit a strong Aharonov-Bohm effect. 
On hyperbolic surfaces with at least two cusps, we show 
that the magnetic Laplacian associated to every magnetic field with 
compact support has purely discrete spectrum for some choices of the vector 
potential, while other choices lead to a situation of limiting absorption
principle. 

We also study perturbations of the metric. 
We show that in the Mourre theory it is not necessary to require 
a decay of the derivatives of the perturbation. This very singular 
perturbation is then brought closer to the perturbation 
of a potential.
\end{abstract}



\section{Introduction} 
Let $X$ be a smooth manifold of dimension $n$, diffeomorphic outside 
a compact set to a cylinder 
$(1,\infty)\times M$, where $M$ is a possibly disconnected closed manifold. 
On $X$ we consider asymptotically conformally 
cylindrical metrics, i.e., perturbations of the
metric given near the border $\{\infty\}\times M$ by: 
\begin{align}\label{mc}
g_p=y^{-2p}(dy^2+h),&&y\to\infty\end{align}
where $h$ is a metric on $M$ and $p>0$. If $p=1$ and $h$ is flat, the ends 
are \emph{cusps}, i.e., complete hyperbolic of finite volume. For $p>1$ one 
gets the (incomplete) metric horns. 
 
The refined properties of the essential spectrum of the Laplace-Beltrami
operator $\Delta_p:=d^*d$ have been studied by Froese and Hislop \cite{FH} 
in the complete case. For the unperturbed metric \eqref{mc}, they get 
\begin{eqnarray}\label{e:Ise}
\sigma_{\rm ess}(\Delta_p)=[\kappa(p), \infty), \mbox{ where } 
\left\{\begin{array}{ll}
\kappa(p)=0, &\mbox{ for } p	\pp 1 \\
\kappa(1)=\left(\frac{n-1}{2}\right)^2. &
\end{array}\right.
\end{eqnarray}
The singular continuous part of the spectrum is empty and the eigenvalues 
distinct from $\kappa(p)$ are of finite multiplicity and may
accumulate only at $\kappa(p)$. Froese and Hislop actually show a 
limiting absorption principle, a stronger result, see also 
\cite{DHS, FHP, FHP1, H} for the continuation of their ideas. 
Their approach relies on a positive commutator technique introduced 
by E.\ Mourre in \cite{mou}, see also \cite{ABG} and references
therein. See for instance \cite{Guillope, Kumura} for different
methods. 

Consider more generally a conformal perturbation of the metric \eqref{mc}. 
Let $\rho\in\cC^\infty(X,\rz)$ be such that $\inf_{y\in X} (\rho(y))\pg -1$. 
Consider the same problem as above for the metric
\begin{equation} \label{e:Ipertu}
\tilde g_p=(1+\rho) g_p, \mbox{ for large } y.
\end{equation}
To measure the size of the perturbation, we compare it to the lengths of 
geodesics. Let $L\in\cC^\infty(X)$ be defined by
\begin{eqnarray}\label{e:IL}
L\geq 1,\quad L(y)=\left\{\begin{array}{ll}
\frac{y^{1-p}}{1-p}&\mbox{ for } p	\pp 1 \\
\ln(y) &\mbox{ for } p	= 1
\end{array}\right., \mbox{ for } y \mbox { big enough}.
\end{eqnarray}

In \cite{FH}, one essentially asks that 
\begin{equation*}
L^2\rho, L^2d\rho \mbox{ and } L^2\Delta_g\rho \mbox{ are in } L^\infty(X).
\end{equation*}
to obtain the absence of singular continuous spectrum and local finiteness 
of the point spectrum. On one hand, one knows from the perturbation of 
a Laplacian by a short-range potential $V$ that only the speed of the
decay of $V$ is important to conserve these properties. On the other
hand, in \cite{GG} and in a general setting, one shows that only the
fact that $\rho$ tends to $0$ is enough to ensure the stability of
the essential spectrum. Therefore, it is natural to ask whether the
decay of the metric (without decay conditions on the derivatives) is
enough to ensure the conservation of these properties. In this paper,
we consider that $\rho=\rho_{\rm sr }+\rho_{\rm lr }$ decomposes in
short-range and long-range components. We ask the long-range component 
to be radial. We also assume that 
there exists $\varepsilon\pg 0$ such that 
\begin{equation}\label{e:Iconf}\begin{split}
L^{1+\varepsilon}\rho_{\rm sr } \text{ and }
d\rho_{\rm sr }, \Delta_g\rho_{\rm sr } \in L^\infty(X),\\
L^{\varepsilon}\rho_{\rm lr }\text{, } 
 L^{1+\varepsilon}d\rho_{\rm lr } \text{ and } \Delta_g \rho_{\rm sr } 
\in L^\infty(X).
\end{split}\end{equation}
Going from $2$ to $1+\varepsilon$ is not a significant 
improvement as it relies on the use of an optimal version of the Mourre 
theory instead of the original theory, see \cite{ABG} and references therein. 
Nevertheless, the fact that the derivatives are asked only 
to be bounded and no longer to decay is a real improvement due to our method. 
We prove this result in Theorem \ref{t:mourre0}. In the Mourre theory, one 
introduces a conjugate operator to study a given operator. The conjugate 
operator introduced in \cite{FH} is too rough to handle very singular 
perturbations. 
In our paper, we introduce a conjugate operator local in energy to avoid 
the problem. We believe that our approach could be implemented easily
in the manifold settings from \cite{bouclet, DHS, FH, FHP, FHP1, H, KT}
to improve results on perturbations of the metric.

A well-known dynamical consequence of the absence of singular continuous 
spectrum and of the local finiteness of the point spectrum is that for 
an interval $\cJ$ that contains no eigenvalue of the Laplacian,
for all $\cchi\in\cC^\infty_c(X)$ and $\phi\in L^2(X)$, the norm
$\|\cchi e^{it\Delta_p} E_\cJ(\Delta_p) \phi\|$ tends to $0$ as $t$ tends to 
$\pm \infty$. In other words, if you let evolve long enough a particle 
which is located at scattering energy, it eventually becomes located very far 
on the exits of the manifold. Add now a magnetic field $B$
with compact support and look how strongly it can interact with the particle. 
Classically there is no interaction as $B$ and the particle 
are located far from each other. One looks for a quantum effect.

The Euclidean intuition tells us that is no essential difference
between the free Laplacian and the magnetic Laplacian $\Delta_A$, where 
$A$ is a magnetic potential arising from a magnetic field $B$ with compact 
support. They still share the spectral properties of absence of 
singular continuous spectrum and local finiteness of the point spectrum, 
although a long-range effect does occur and destroys the asymptotic 
completeness of the couple $(\Delta, \Delta_A)$; one needs to modify 
the wave operators to compare the two operators, see
\cite{LT}. However, we point out in this paper that the situation
is dramatically different in particular on hyperbolic manifolds
of finite volume, even if the magnetic field is very small in size
and with compact support. We now go into definitions and describe our results.

A magnetic field $B$ is a smooth real
exact $2$-form on $X$. There exists a real $1$-form $A$, called vector 
potential, satisfying $dA=B$. Set $\dA:=d+iA\wedge: \cunc(X)\to\cunc(X,T^*X)$.
The magnetic Laplacian on $\cunc(X)$ is given by $\Delta_A:=\dA^*\dA$. 
When the manifold is complete, $\Delta_A$ is known to be essentially
self-adjoint, see \cite{shubin}. Given two vector potentials $A$ and $A'$ such that $A-A'$ is 
exact, the two magnetic Laplacians $\Delta_A$ and $\Delta_{A'}$ are 
unitarily equivalent, by gauge invariance. Hence when $H^1_\dR(X)=0$, 
the spectral properties of the magnetic Laplacian do not depend on the choice 
of the vector potential, so we may write $\Delta_B$ instead of 
$\Delta_A$. 

The aim of this paper is the study of the spectrum of magnetic 
Laplacians on a manifold $X$ with the metric \eqref{pme}, which includes
the particular case \eqref{mc}. In this introduction we restrict 
the discussion to the complete case, i.e.\ $p\leq 1$. We focus first on 
the case of gauge invariance, i.e., $H^1_\dR(X)=0$, and we simplify 
the presentation assuming that the boundary is connected. 
We classify magnetic fields.

\begin{definition}\label{d:intro}
Let $X$ be the interior of a compact manifold 
with boundary $\oX$. Suppose that $H^1_\dR(X)=0$ and that 
$M=\partial \oX$ is connected. Let $B$ be a magnetic field on $X$
which extends smoothly to a $2$-form on $\oX$. 
We say that $B$ is \emph{trapping} if 
\begin{itemize} 
\item[1)] either $B$ does not vanish identically on $M$, or
\item[2)] $B$ vanishes on $M$ but defines a non-integral cohomology class 
 $[2\pi B]$ inside the relative cohomology group $H^2_\dR(X,M)$. 
\end{itemize}
\end{definition}
Otherwise, we say that $B$ is \emph{non-trapping}. 
This terminology is motivated by the spectral consequences a) and c)
of Theorem \ref{t:Ithmag}.
The definition can be generalized to the case where $M$ is
disconnected (Definition \ref{magtrap}). The condition of $B$ being
trapping can be expressed in terms of any vector potential $A$ (see
Definition \ref{def7} and Lemma \ref{l:compsupp}). When
$H^1_\dR(X)\neq 0$, the trapping condition makes sense only for
vector potentials, see Section \ref{s:mfc} and Theorem \ref{t:IAB}.

Let us fix some notation. Given two Hilbert spaces $\Hr$ and $\Kr$, we
denote by $\Bc(\Hr, \Kr)$ and $\Kc(\Hr, \Kr)$ the bounded and compact
operators acting from $\Hr$ to $\Kr$, respectively. Given $s\geq 0$,
let $\Lr_{s}$ be the domain of $L^s$ equipped with the graph norm. We
set $\Lr_{-s}:=\Lr_{s}^*$ where the adjoint space is defined so that 
$\Lr_{s}\subset L^2(X, g_p) \subset \Lr_{s}^*$, using the
Riesz lemma. Given a subset $I$ of $\R$, let $I_{\pm}$ be the set of
complex numbers $x\pm iy$, where $x\in I$ and $y>0$. For simplicity, 
in this introduction we state our result only for the unperturbed 
metric \eqref{mc}.

\begin{theorem}\label{t:Ithmag}
Let $1\geq p>0$, $g_p$ the metric given by \eqref{mc}. 
Suppose that $H_1(X,\zz)=0$ and that $M$ is connected. 
Let $B$ be a magnetic field which extends smoothly to $\oX$.
If $B$ is trapping then:
\begin{itemize} 
\item[a)] The spectrum of $\Delta_B$ is purely discrete.
\item[b)] The asymptotic of its eigenvalues is given by 
\begin{equation}\label{e:Ithmag}
N_{B,p}(\lambda) \approx \begin{cases}
C_1\lambda^{n/2}&
\text{for $1/n< p$,}\\
C_2\lambda^{n/2}\log \lambda &\text{for $p=1/n$,}\\
C_3\lambda^{1/2p}&\text{for $0<p<1/n$}
\end{cases}
\end{equation}
in the limit $\lambda\to\infty$, where $C_3$ is given in Theorem \ref{t:thmag}, and
\begin{equation}\label{e:Ic}
C_1=\frac{\Vol(X,g_p)\Vol(S^{n-1})}{n(2\pi)^n},\quad 
C_2=\frac{\Vol(M,h)\Vol(S^{n-1})}{2(2\pi)^n}. 
\end{equation} 
\end{itemize} 
If $B$ is non-trapping with compact support in $X$ then 
\begin{itemize} 
\item[c)] The essential spectrum of $\Delta_B$ is $[\kappa(p), \infty)$.
\item[d)] The singular continuous spectrum of $\Delta_B$ is empty. 
\item[e)] The eigenvalues of $\Delta_B$ are of finite multiplicity and can accumulate only in $\{\kappa(p)\}$.
\item[f)] Let $\cJ$ a compact interval such that $\cJ\cap \big(\{\kappa(p)\}\cup
 \sigma_{\rm pp}(H)\big)=\emptyset$. Then, for all $s\in ]1/2, 3/2[$ and all $A$ such that
 $dA=B$, there is $c$ such that 
\begin{equation*}
\|(\Delta_A-z_1)^{-1} - (\Delta_A-z_2)^{-1} \|_{\Bc(\Lr_s, \Lr_{-s})} 
\leq c \|z_1-z_2\|^{s-1/2},
\end{equation*}
for all $z_1, z_2 \in \cJ_{\pm}$.
\end{itemize} 
\end{theorem}
The statements a) and b) follow from 
general results from \cite{wlom}. This part relies on the Melrose 
calculus of cusp pseudodifferential operators (see e.g., \cite{meni96c})
and is proved in Theorem \ref{t:thmag} for the perturbed metric \eqref{pme}.
We start from the basic observation that for smooth vector potentials, 
the magnetic Laplacian belongs to the cusp calculus with positive weights. 
For this part, we can treat the metric \eqref{pme} which is
quasi-isometric (but not necessarily asymptotically equivalent) to \eqref{mc}. 
Moreover, the finite multiplicity of the point spectrum 
(which is possibly not locally finite) in e) follows from Appendix
\ref{fm} for this class of metrics. 

The statement $c)$ follows directly from the analysis of
the free case in Section \ref{s:free}. The perturbation of the metric 
is considered in Proposition \ref{p:thema} and relies on general results 
on stability of the essential spectrum shown in \cite{GG}. 
The points d), f) and e) rely on the use of an optimal version of 
Mourre theory, see \cite{ABG}. They are developed in 
Theorem \ref{t:mourre0} for perturbations satisfying \eqref{e:Iconf}. 
Scattering theory under short-range perturbation of a potential and of 
a magnetic field is also considered.

The condition of being trapping (resp.\ non-trapping) is discussed 
in Section \ref{s:inv} and is equivalent to having empty 
(resp.\ non-empty) essential spectrum in the complete case. 
The terminology arises from the dynamical consequences of this theorem
and should not be confused with the classical terminology. 
Indeed, when $B$ is trapping, the spectrum of $\Delta_B$ is purely discrete 
and for all non-zero $\phi$ in $L^2(X)$, there is $\cchi\in\cC^\infty_c(X)$ 
such that $1/T \int_{0}^T \|\cchi e^{it \Delta_B} \phi\|^2 dt$ tends to a 
non-zero constant as $T$ tends to $\pm \infty$. On the other hand, taking 
$\cJ$ as in f), for all $\cchi\in\cC^\infty_c(X)$ one gets that
$\cchi e^{it \Delta_B}E_\cJ(\Delta_A) \phi$ tends to zero, when $B$ 
is non-trapping and with compact support. 

If $H^1_\dR(M)\neq 0$ (take $M=S^1$ for instance), there exist 
some trapping magnetic fields with \emph{compact support}. 
We construct an explicit example in Proposition \ref{p:contreex}. 
We are able to construct some examples in dimension $2$ and higher than $4$ 
but there are topological obstructions in dimension $3$, see section 
\ref{s:inv}. As pointed out above regarding the Euclidean case, the fact 
that a magnetic field with compact support can turn off the essential spectrum 
and even a situation of limiting absorption principle is somehow
unexpected and should be understood as a strong long-range effect.

We discuss other interesting phenomena in Section 
\ref{s:coupling}. Consider $M=S^1$ and take a trapping magnetic field 
$B$ with compact support and a coupling constant $g\in\R$. Now remark 
that $\Delta_{g B}$ is non-trapping if and only if $g$ belongs to the
discrete group $c_B\Z$, for 
a certain $c_B\neq 0$. When $g\notin c_B\Z$ and $p\geq 1/n$, 
the spectrum of $\Delta_{g B}$ is discrete and the eigenvalue asymptotics 
do not depend either on $B$ or on $g$. It would be very interesting
to know whether the asymptotics of embedded eigenvalues, or more likely
of resonances, remain the same when $g\in c_B\Z$, 
(see \cite{tanya} for the case $g=0$). It would be also interesting to
study the inverse spectral problem and ask if the magnetic field
could be recovered from the knowledge of the whole spectrum, since
the first term in the asymptotics of eigenvalues does not feel it. 
 
Assume now that gauge invariance does not hold, i.e., $H^1_\dR(X)\neq 0$. 
In quantum mechanics, it is known that the choice of a vector
potential has a physical meaning. This is known as the Aharonov-Bohm
effect \cite{AB}. Two choices of magnetic potential may
lead to in-equivalent magnetic Laplacians. In $\rz^2$ with a 
bounded obstacle, this phenomenon can be seen through a difference 
of wave phase arising from two
non-homotopic paths that circumvent the obstacle. Some long-range effect 
appears, for instance in the scattering matrix like in \cite{R, RY, RY2}, 
in an inverse-scattering problem \cite{N,W} or in the semi-classical regime 
\cite{BR}. See also \cite{helffer} for the influence of the obstacle on the bottom of the spectrum. In all of these cases, the essential spectrum remains the same.

In section \ref{s:AB}, we discuss the Aharonov-Bohm effect in our setting. 
In light of Theorem \ref{t:Ithmag}, one expects a drastic effect. 
We show that the choice of a vector potential can indeed have a significant 
spectral consequence. For one choice of vector potential, the essential spectrum 
could be empty and for another choice it could be a half-line. 
This phenomenon is generic for hyperbolic surfaces of finite 
volume, and also appears for hyperbolic $3$-manifolds. 
We focus the presentation on magnetic fields $B$ with compact support. 
We say that a smooth vector potential $A$ (i.e., a smooth $1$-form 
on $\oX$) is trapping if $\Delta_{A}$ has compact resolvent, and
non-trapping otherwise. By Theorems \ref{t:thmag} and 
\ref{t:mourre0}, for $p\leq 1$, this is equivalent to 
Definition \ref{def7}. It also follows that when the metric is of 
type \eqref{mc}, $A$ is trapping if 
and only if a)--b) of Theorem \ref{t:Ithmag} hold for $\Delta_{A}$, while 
$A$ is non-trapping if and only if $\Delta_{A}$ satisfies
c)--f) of Theorem \ref{t:Ithmag}.

\begin{theorem}\label{t:IAB}
Let $X$ be a complete oriented hyperbolic surface of finite volume and $B$ a 
smooth magnetic field on the compactification $\oX$.
\begin{itemize}
\item If $X$ has at least $2$ cusps, then for all
$B$ there exists both trapping and non-trapping vector 
potentials $A$ such that $B=dA$.
\item If $X$ has precisely $1$ cusp, choose $B=dA=dA'$ where $A,A'$ are 
smooth vector potentials for $B$ on $\oX$. 
Then 
\begin{align*}
\text{$A$ is trapping $\Longleftrightarrow A'$ 
is trapping $\Longleftrightarrow \int_X B\in 2\pi \zz$.}
\end{align*}
\end{itemize}
\end{theorem}
This follows from Corollary \ref{81}. 
More general statements are valid also in dimension $3$, see 
Section \ref{hyp}.

This implies on one hand that for a choice of $A$, as one has the 
points c), d) and e), a particle located at a scattering energy
escapes from any compact set; on the other hand taking a trapping
choice, the particle will behave like an eigenfunction and will
remain bounded. It is interesting that the dimension $3$ is
exceptional in the Euclidean case \cite{Y} and that we are able to
construct examples of such a behavior in any dimension. 

In the first appendix, we discuss the key notion of $C^1$ regularity
for the Mourre theory and make it suitable to the manifold context and 
for our choice of conjugate operator. As pointed in \cite{GG0}, this is a 
key hypothesis in the Mourre theory in order to apply the Virial theorem 
and deduce the local finiteness of the embedded eigenvalues. 
In the second appendix, we recall that (cusp) elliptic, not necessarily 
fully elliptic, cusp operators have $L^2$ eigenvalues of finite multiplicity. 
Finally, in the third appendix, we give a criteria of stability of the 
essential spectrum, by cutting a part of the space, encompassing
incomplete manifolds. Some of the results concerning the essential
spectrum and the asymptotic of eigenvalues were already present in the
unpublished preprint \cite{GMo}.

\subsection*{Acknowledgements.}
We acknowledge useful discussions with Barbu Berceanu,
Dan Burghelea, Jan De\-re\-zi\'nski, Vladimir Georgescu,
Bernard Helffer, Andreas Knauf, Fran\c cois Ni\-co\-leau, Marius M\u antoiu and
Radu Purice. We are also grateful to the referee for helpful remarks. 

The authors were partially supported from the contract MERG 006375,
funded by the European Commission. The second author was partially 
supported by the contracts 2-CEx06-11-18/2006 and CNCSIS-GR202/19.09.2006 
(Romania).

\section{Cusp geometry}\label{lcm}
\subsection{Definitions}\label{s:gene}
This section follows closely \cite{wlom}, see also \cite{meni96c}.
Let $\oX$ be a smooth $n$-dimensional compact manifold
with closed boundary $M$, and $x:\oX\to[0,\infty)$ a
boundary-defining function. A \emph{cusp metric} on $\oX$ is a
complete Riemannian metric $g_0$ on $X:=\oX\setminus M$
which in local coordinates near the boundary takes the form
\begin{equation}\label{cume}
g_0=a_{00}(x,y)\frac{dx^2}{x^4}+\sum_{j=1}^{n-1} a_{0j}(x,y)
\frac{dx}{x^2}dy_j+\sum_{i,j=1}^{n-1}a_{ij}(x,y)dy_idy_j
\end{equation}
such that the matrix $(a_{\alpha\beta})$ is smooth and non-degenerate
down to $x=0$. For example, if $a_{00}=1$, $a_{0j}=0$ and
$a_{ij}$ is independent of $x$, we get a product metric near $M$.
If we set $y=1/x$, a cusp metric is nothing but 
a quasi-isometric deformation of a
cylindrical metric, with an asymptotic expansion for 
the coefficients in powers of $y^{-1}$.
We will focus on the \emph{conformally cusp} metric
\begin{align}\label{pme}
g_p:=x^{2p}g_0,\end{align}
where $p>0$. Note that \eqref{mc} is a particular case of such metric. 

Let $\cI\subset\cun(\oX)$ be the principal ideal generated by the function
$x$. Recall \cite{meni96c} that a \emph{cusp vector field} is a smooth 
vector field $V$ on $\oX$ such that $dx(V)\in\cI^2$. 
The space of cusp vector fields forms a
Lie subalgebra $\ccV$ of the Lie algebra $\cV$ of smooth vector fields on $\oX$. 
In fact, there exists a natural vector bundle $\ctX$ over $\oX$ whose space of 
smooth sections is $\ccV$, and a natural map $\ctX\to TX$ which induces 
the inclusion $\ccV\hookrightarrow \cV$.
Let $E,F\to\oX$ be smooth vector bundles. The space of cusp differential
operators $\Diffc(\oX,E,F)$ is the space of those differential operators
which in local trivializations can be written as composition of cusp
vector fields and smooth bundle morphisms down to $x=0$.

The \emph{normal operator} of $P\in\Diffc(\oX,E,F)$ is the family of operators
defined by
\[\rz\ni\xi\mapsto\cN(P)(\xi):=
\left(e^{i\xi/x}Pe^{-i\xi/x}\right)_{|x=0}\in\Diff(M,E_{|M},F_{|M}).\] 
\begin{example}\label{ex1}
$\cN(x^2\px)(\xi)=i\xi$.
\end{example} 
Note that $\ker \cN=\cI\cdot\Diffc$, which we denote again by 
$\cI$. The normal operator map is linear and multiplicative.
It is also invariant under the conjugation by
powers of $x$. Namely, if $P\in\Diffc$ and $s\in\cz$ then
$x^s Px^{-s}\in\Diffc$ and $\cN(x^s Px^{-s})=\cN(P)$. Concerning 
taking the (formal) adjoint, one needs to specify the volume form 
on the boundary. 

\begin{lemma} \label{lema6}
Let $P\in\Diffc(\oX,E,F)$ be a cusp operator and $P^*$ its adjoint with
respect to $g_0$. Then $\cN(P^*)(\xi)$ is the adjoint of $\cN(P)(\xi)$
with respect to the metric on $E_{|M},F_{|M}$ induced by restriction,
for the volume form ${a_{0}}^{1/2} \vol_{h_0}$, where the metric $h_0$ on $M$
is defined from rewriting $g_0$ as in \eqref{cume2}. 
\end{lemma}
The principal symbol of a cusp operator on $X$ extends as a map on the cusp
cotangent bundle down to $x=0$. This implies that a cusp operator of positive
order cannot be elliptic at $x=0$ in the usual sense. A cusp operator 
is called cusp-elliptic if its
principal symbol is invertible on $\ctsX\setminus\{0\}$ down to $x=0$.

\begin{definition} 
A cusp operator is called \emph{fully elliptic} if it is cusp-elliptic and 
if its normal operator is invertible for all values of $\xi\in\rz$.
\end{definition} 
An operator $H\in x^{-l}\Diffc^k(\oX,E,F)$ is called 
a cusp differential operator of type 
$(k,l)$.

Fix a product decomposition of $X$ near $M$, compatible with the 
boundary-defining function $x$. This gives a splitting of the cusp 
cotangent bundle on $\oX$ in a neighborhood of $M$: 
\begin{equation}\label{decom}
\ctsX\simeq T^*M\oplus\langle x^{-2}dx\rangle.
\end{equation}

\begin{lemma}\label{prop3}
The de Rham differential $d:\cun(X)\to \cun(X,T^*X)$ restricts to a 
cusp differential operator $d:\cun(\oX)\to\cun(\oX,\ctsX)$. 
Its normal operator in the decomposition \eqref{decom} is
\[\cN(d)(\xi)=
\begin{bmatrix}d^M\\i\xi\end{bmatrix}\]
where $d^M$ is the partial de Rham differential in the $M$ factor of the
product decomposition.
\end{lemma}
\proof 
Let $\omega\in\cun(\oX)$ and decompose $d\omega$ according to \eqref{decom}:
\[d\omega=d^M\omega + \px(\omega) dx=d^M\omega+x^2\px(\omega)\frac{dx}{x^2}.\]
Since $d^M$ commutes with $x$, it follows from the definition that
$\cN(d^M)=d^M$. The result follows using Example \ref{ex1}.
\qed

\subsection{Relative de Rham cohomology}\label{rdR}

Recall \cite{BT} that the cohomology of $\oX$ and the relative 
cohomology groups of 
$(\oX,M)$ (with real coefficients) can be computed using smooth differential
forms as follows let $\Lambda^*(\oX)$ denote the space of forms
smooth on $\oX$ down to the boundary. Let $\Lambda^*(\oX,M)$ denote
the subspace of those forms whose pull-back to $M$ vanishes. These
spaces form complexes for the de Rham differential (because $d$
commutes with pull-back to $M$) and their quotient is the de Rham
complex of $M$: 
\[0\to \Lambda^*(\oX,M)\hookrightarrow
\Lambda^*(\oX)\to\Lambda^*(M)\to 0.\] 
The induced long exact sequence in cohomology 
is just the long exact sequence of the pair $(\oX,M)$. 

\subsection{Cusp de Rham cohomology}

Notice that $d$ preserves the space of cusp differential forms. 
Indeed, since $d$ is a derivation and using Lemma \ref{prop3},
it suffices to check this property for a set of local generators of 
$\cun(\oX,\ctsX)$. Choose local coordinates $(y_j)$ on $M$
and take as generators $x^{-2}dx$ and $dy_j$, which are closed.

Let $\cH^*(X)$ denote the cohomology of the complex of cusp differential
forms $(\cun(\oX,\Lambda^*(\ctX)), d)$ with respect to the de Rham differential.

\begin{proposition}\label{dec}
$\cH^k(X)=H^k(X)\oplus H^{k-1}(M)^2$.
\end{proposition}
\proof 
The short exact sequence of de Rham complexes
\begin{equation*} 
0\to \Lambda^*(T\oX)\hookrightarrow 
\Lambda^*(\ctX)\to \Lambda^{*-1}(TM)^2\to 0,
\end{equation*} 
where the second map is given by
\begin{equation}\label{restfr}
\Lambda^k(\ctX)\ni\omega\mapsto \left((x^2\px\lrcorner \omega)_{x=0},
(\partial_x(x^2\px\lrcorner \omega))_{x=0}\right),
\end{equation}
gives rise to a long exact sequence in cohomology. Now the composition
\[\Lambda^*(T\oX)\hookrightarrow \Lambda^*(\ctX)\hookrightarrow \Lambda^*(TX)\]
is a quasi-isomorphism, since de Rham cohomology can be computed either with 
smooth forms on $\oX$, or with smooth forms on $X$. Thus in cohomology 
the map induced from $\Lambda^*(T\oX)\hookrightarrow 
\Lambda^*(\ctX)$ is injective.
\qed

\section{The magnetic Laplacian}\label{sectionmagn}
\subsection{The magnetic Laplacian on a Riemannian manifold}
A \emph{magnetic field} $B$ on the Riemannian manifold $(X,g)$
is an exact real-valued $2$-form. A \emph{vector potential}
$A$ associated to $B$ is a $1$-form such that $dA=B$. We
form the \emph{magnetic Laplacian} acting on $\cun(X)$:
\[\Delta_A:=\dA^*\dA.\]
This formula makes sense for complex-valued $1$-forms $A$.
Note that when $A$ is real, $\dA$ is a metric connection on the trivial
bundle $\underline{\cz}$ with the canonical metric, and $\Delta_A$
is the connection Laplacian.

If we alter $A$ by adding to it a real exact form,
say $A'=A+df$, the resulting magnetic Laplacian satisfies
\[\Delta_{A'}=e^{-if}\Delta_Ae^{if}\]
so it is unitarily equivalent to $\Delta_A$ in $L^2(X,g)$.
Therefore if $H^1_{\mathrm{dR}}(X)=0$ (for instance if $\pi_1(X)$ is
finite; see \cite{BT}) then $\Delta_A$ depends, up to unitary
equivalence, only on the magnetic field $B$. This property is called
\emph{gauge invariance}. For a more refined analysis of gauge
invariance, see \cite{gruber}. 

One usually encounters gauge invariance as a
consequence of $1$-connectedness (i.e., $\pi_1=0$). But in dimensions
at least $4$, every finitely presented group (in particular, every
finite group) can be realized as $\pi_1$ of a compact manifold. Thus
the hypothesis $\pi_1=0$ is unnecessarily restrictive, it is enough to assume
that its abelianisation is finite. 

While the properties of $\Delta_A$ in $\rz^n$ with the flat metric
are quite well understood, the (absence of) essential spectrum of
magnetic Laplacians on other manifolds has not been much studied so far.
One exception is the case of bounded geometry, studied in \cite{KS}.
However our manifolds are not of bounded geometry because the injectivity
radius tends to $0$ at infinity. 
 
\subsection{Magnetic fields and cohomology}\label{s:mfc} 
Recall that $A$ is a (smooth) \emph{cusp $1$-form} (not to be confused 
with the notion of cusp form from automorphic form theory) on $\oX$ if 
$A\in\cun(X,T^*X)$ is a real-valued $1$-form satisfying near $\partial X$
\begin{equation}\label{aft}
A=\varphi(x)\frac{dx}{x^2}+\theta(x) 
\end{equation}
where $\varphi\in\cun(\oX)$ and $\theta\in\cun([0,\varepsilon)\times M,
\Lambda^1(M))$, or equivalently $A$ is a smooth section in $\ctsX$
over $\oX$.

\begin{proposition}\label{p:cuspduB} 
Let $B$ be a cusp $2$-form. Suppose that $B$ is exact on $X$, and its image
by the map \eqref{restfr} is exact on $M$. Then there exists
a smooth cusp $1$-form $A$ on $\oX$ such that $dA=B$.
\end{proposition} 
\proof
Note that $B$ is exact as a form on $X$, so $dB=0$ on $X$.
By continuity, $dB=0$ on $\oX$ (in the sense of cusp forms)
so $B$ defines a cusp cohomology $2$-class. 
By hypothesis, this class maps to $0$ by restriction to $X$.
Now the pull-back of $B$ to the level surfaces $\{x=\varepsilon\}$ is closed;
by continuity, the image of $B$ through the map \eqref{restfr} is closed on 
$M$. Assuming that this image is exact, it follows from Proposition
\ref{dec} that $B$ is exact as a cusp form.
\qed
 
By Lemmata \ref{lema6} and \ref{prop3},
$\Delta_A$ is a cusp differential operator of order $(2p,2)$. 
 
\begin{definition} \label{def7}
Let $A$ be a (complex-valued) cusp vector potential. 
Given a connected component $M_0$ of $M$, 
we say $A$ is a \emph{trapping} vector potential on $M_0$ if
\begin{itemize} 
\item either the restriction $\varphi_0:=\varphi(0)$ is not 
constant on $M_0$, 
\item or $\theta_0:=\theta(0)$ is not closed on $M_0$, 
\item or the cohomology class
$[{\theta_0}_{|M_0}]\in H^1_{\mathrm{dR}}(M_0)$ does not belong to the image of
\[2\pi H^1(M_0;\zz)\to H^1(M_0,\cz)\simeq H^1_{\mathrm{dR}}(M_0)\otimes\cz.\] 
\end{itemize}
and \emph{non-trapping on $M_0$} otherwise.
\end{definition} 

We say that $A$ is \emph{trapping} if it is trapping on each 
connected component of $M$. The vector potential is said to be 
\emph{non-trapping} if it is non-trapping on at least one connected 
component of $M$. 
If $A$ is non-trapping on all connected component of $M$, we say 
that it is \emph{maximal non-trapping}.

\begin{remark}
The trapping notion can be expressed solely in terms of the magnetic field $B
=dA$ when $H^1(X)=0$, see Lemma \ref{l:compsupp}.
\end{remark}

We comment briefly the terminology. When $A$ is constant in $x$
near $M$, the multiplicity of the absolutely continuous part of the spectrum 
of $\Delta_A$ will be given by the number of connected component of $M$ 
on which $A$ is non-trapping. Hence, taking $A$ maximal non-trapping 
maximizes the multiplicity of this part of the spectrum.

We refer to \cite{BT} for an exposition of cohomology with integer
coefficients. The trapping property is determined only by the
asymptotic behavior of $A$. More precisely, if $A'$ is also of the
form \eqref{aft} with $\varphi(0)=0$ and $\theta(0)=0$ then $A$ is a
(non-)trapping vector potential if and only if $A+A'$ is. 

The term ``trapping" is motivated by dynamical consequences
of Theorems \ref{p:pmag} and \ref{t:mourre0} and has nothing
to do with the classical trapping condition. This terminology is
also supported by the examples given in Section \ref{s:nonstab}.

For a trapping vector potential, $x^{2p}\Delta_A$ is a
fully-elliptic cusp operator. In turn, this implies that $\Delta_A$ has 
empty essential spectrum so from a dynamical point of view, a 
particle can not diffuse, in other words it is trapped in the interior of
$X$. Indeed, given a state $\phi\in L^2(X)$, there exists $\cchi$
(the characteristic function of a compact subset of $X$)
such that $1/T\int_0^T \|\cchi e^{it\Delta_A} \phi\|^2 dt$ tends to a positive
constant as $T$ goes to infinity.

On the other hand, if $A$ is a non-trapping vector potential, 
then $\Delta_A$ is not Fredholm between the appropriate cusp Sobolev spaces. 
If the metric is an exact cusp metric and complete, we show that $\Delta_A$
has nonempty essential spectrum also as an unbounded operator in
$L^2$, given by $[\kappa(p), \infty)$ by Proposition \ref{p:thema}. We
go even further and under some condition of decay of $\varphi$ and
$\theta$ at infinity, we show that there is no singular continuous
spectrum for the magnetic Laplacian and that the eigenvalues of
$\R\setminus \{\kappa(p)\}$ are of finite multiplicity and can accumulate
only in $\{\kappa(p)\}$. Therefore given a state $\phi$ which is not an
eigenvalue of $\Delta_A$, one obtains that for all
$\cchi$, $\cchi e^{it\Delta_A}\phi$ tends to $0$ as $t\to\infty$.

When $M$ is connected, the class of non-trapping vector potentials is a group 
under addition but that of trapping vector potential is not. 
When $M$ is disconnected, none of these classes is closed under addition.
Directly from the definition, we get however:

\begin{remark}\label{addtrap} 
Let $A$ be a maximal non-trapping vector potential and let $A'$ be 
a $1$-form smooth up to the boundary. 
Then $A'$ is trapping if and only if $A+A'$ is.
\end{remark}

Let $B$ be a smooth magnetic field on $\oX$ (i.e., a $2$-form)
whose pull-back to $M$ vanishes. 
Since $B$ is exact, it is also closed, thus it defines a relative de
Rham class as in Subsection \ref{rdR}. If this class vanishes, we
claim that there exists a vector potential $A$ for $B$ which is
maximal non-trapping. Indeed, let $A\in \Lambda^1(\oX,M)$ be any (relative)
primitive of $B$. Then $A$ is clearly a cusp form, the singular term
$\phi(0)$ vanishes, and the pull-back of $A$ to each boundary component
vanishes by definition, in particular it defines the null $1$-cohomology
class. From Remark \ref{addtrap} we get 

\begin{corollary}\label{cohcl}
Let $B$ be a cusp magnetic field. Let $B'$ be a smooth magnetic field on 
$\oX$ which vanishes on the boundary and 
which defines the zero relative cohomology class in $H^2(\oX,M)$. Then
$B$ admits (non-)trapping vector potentials if and only if $B+B'$ does.
\end{corollary}

Note that when $H^1(\oX)\neq 0$, a given magnetic field may admit 
both trapping and non-trapping vector potentials. See Theorem \ref{t:IAB} and 
Sections \ref{s:AB} and \ref{hyp}.

\section{The trapping case}
\subsection{The absence of essential spectrum}
In this section, given a smooth cusp $1$-form, we discuss the link
between its behavior at infinity and its trapping properties. 
 
\begin{theorem}\label{p:pmag}
Let $p>0$, $g_p$ a metric on $X$ given by \eqref{pme} near $\partial X$
and $A$ a smooth cusp $1$-form given by \eqref{aft}. Then $\Delta_A$ is 
a weighted cusp differential 
operator of order $(2p,2)$. If $A$ is trapping then 
$\Delta_A$ is essentially self-adjoint on $\cC_c^\infty(X)$, 
it has purely discrete spectrum and its domain is $x^{2p}H^2(X,g_p)$.
\end{theorem}
If $p\leq 1$ then $g_p$ is complete so $\Delta_A$ is essentially
self-adjoint \cite{shubin}. This fact remains true for a trapping $A$
in the incomplete case, i.e.\ $p>1$.

\proof
Using Lemma \ref{prop3}, we get
\[\cN(\dA)(\xi)=
\begin{bmatrix} d^M+i\theta_0\\
i(\xi+\varphi_0)\end{bmatrix}.\]
Suppose that $x^{2p}\Delta_A$ is not fully elliptic, so there exists 
$\xi\in\rz$ and $0\neq u\in\ker(\cN(x^{2p}\Delta_A)(\xi))$. By elliptic 
regularity, $u$ is smooth.
We replace $M$ by one of its connected components on which $u$ does not
vanish identically, so we can suppose that $M$ is connected.
Using Lemma \ref{lema6}, by integration by parts
with respect to the volume form ${a_0}^{1/2}dh_0$ on $M$
and the metric $h_0$ on $\Lambda^1(M)$, we see
that $u\in\ker(\cN(\Delta_A)(\xi))$ implies $u\in\ker(\cN(\dA)(\xi))$.
Then
\begin{align}\label{xvp}
(\xi+\varphi_0)u= 0\mbox{ and }(d^M+i\theta_0)u= 0, 
\end{align}
so $u$ is a global parallel section in the trivial bundle $\ucz$
over $M$, with respect to the connection $d^M+i\theta_0$.
This implies
\begin{equation*}\label{zudg}
0=(d^M)^2 u=d^M(-iu\theta_0)
=-i(d^M u)\wedge\theta_0-iud^M\theta_0=-iud^M\theta_0.
\end{equation*} 
By uniqueness of
solutions of ordinary differential equations, $u$ is never $0$, so
$d^M\theta_0=0$. Furthermore, from \eqref{xvp}, we see that $\varphi_0$ 
equals the constant function $-\xi$. It remains to 
prove the assertion about the cohomology class $[\theta_0]$.

Let $\tM$ be the universal cover of $M$. Denote by ${\tilde{u}}$, ${\tilde{\theta}_0}$ the
lifts of $u,\theta_0$ to $\tM$. The equation $(d^M+i\theta_0)u=0$ lifts to
\begin{equation}\label{ecu}
(d^{\tM}+i{\tilde{\theta}_0}){\tilde{u}}=0.
\end{equation}
The $1$-form ${\tilde{\theta}_0}$ is closed on the simply connected manifold
$\tM$, hence it is exact (by the universal coefficients formula,
$H^1(\tM,\cz)=H_1(\tM,\cz)=H_1(\tM;\zz)\otimes\cz$, and $H_1(\tM;\zz)$ 
vanishes as it is the abelianisation of $\pi_1(\tM)$). Let $v\in\cun(\tM)$
be a primitive of ${\tilde{\theta}_0}$, i.e., 
$d^{\tM}v={\tilde{\theta}_0}$. Then, from \eqref{ecu}, ${\tilde{u}}=Ce^{-iv}$ 
for some constant $C\neq 0$.

The fundamental group $\pi_1(M)$ acts to the right on $\tM$ via
deck transformations. The condition that ${\tilde{u}}$ be the lift of
$u$ from $M$ is the invariance under the action of $\pi_1(M)$, in other words
\[{\tilde{u}}(y)={\tilde{u}}(y[\gamma])\]
for all closed loops $\gamma$ in $M$. This is obviously equivalent to
\begin{align*}
v(y[\gamma])-v(y)\in 2\pi\zz,&&\forall y\in \tilde{M}.\end{align*}
Let $\tilde{\gamma}$ be the lift of $\gamma$ starting in $y$. Then
\begin{align*}
v(y[\gamma])-v(y)=&\int_{\tilde{\gamma}}d^{\tM}v
=\int_{\tilde{\gamma}}{\tilde{\theta}_0}=\int_\gamma \theta_0.
\end{align*}
Thus the solution ${\tilde{u}}$ is $\pi_1(M)$-invariant if and only if the 
cocycle $\theta_0$ evaluates to an integer multiple of $2\pi$ on each closed 
loop $\gamma$. These loops span $H_1(M;\zz)$, so $[\theta_0]$ lives in the 
image of $H^1(M;\zz)$ inside $H^1(M,\cz)=\Hom(H_1(M;\zz),\cz)$. Therefore the 
solution $u$ must be identically $0$ unless $\varphi_0$ is constant, 
$\theta_0$ is closed and $[\theta_0]\in 2\pi H^1(M;\zz)$.

Conversely, if $\varphi_0$ is constant, $\theta_0$ is closed and
$[\theta_0]\in 2\pi H^1(M;\zz)$ then ${\tilde{u}}=e^{-iv}$ as above is 
$\pi_1(M)$-invariant, so it is the lift of some $u\in\cun(M)$ which
belongs to $\ker(\cN(x^{2p}\Delta_A)(\xi))$ for $-\xi$ equal to the
constant value of $\varphi_0$. 

The conclusion of the theorem is now a consequence of general properties of 
the cusp calculus \cite[Theorem 17]{wlom}. Namely, since 
$\Delta_A$ is fully elliptic, there 
exists an inverse in $x^{2p}\Psi_c^{-2}(X)$ (a micro-localized version of 
$\Diff_c(X)$) modulo compact operators. If $p>0$, this pseudo-inverse is 
itself compact. The operators in the cusp calculus act 
by closure on a scale of Sobolev spaces. It follows easily that for
$p\geq 0$, a symmetric fully elliptic cusp operator in
$x^{2p}\Psi_c^{2}(X)$ is essentially self-adjoint, with domain
$x^{2p}H^2_c(X)$. Thus $\Delta_A$ is self adjoint with compact inverse
modulo compact operators, which shows that the spectrum is purely discrete.
\qed

The cusp calculus \cite{meni96c} is a particular instance of Melrose's program of 
micro-localizing boundary fibration structures. It is a special case of the 
fibered-cusp calculus \cite{mame99} and can be obtained using
the groupoid techniques of \cite{lani}.

\subsection{Eigenvalue asymptotics for trapping magnetic Laplacians}

If $A$ is a trapping vector potential, the associated magnetic Laplacian has
purely discrete spectrum. In this case we can give the first term
in the eigenvalue growth law.

\begin{theorem}\label{t:thmag}
Let $p>0$, $g_p$ a metric on $X$ given by \eqref{pme} near $\partial X$
and $A\in\cun(X,T^*X)$ a 
complex-valued 
trapping vector potential in the sense of Definition \ref{def7}.
Then the eigenvalue counting function of $\Delta_A$ satisfies
\begin{equation}\label{e:thmag}
N_{A,p}(\lambda) \approx \begin{cases}
C_1\lambda^{n/2}&
\text{for $1/n< p<\infty$,}\\
C_2\lambda^{n/2}\log \lambda &\text{for $p=1/n$,}\\
C_3\lambda^{1/2p}&\text{for $0<p<1/n$}
\end{cases}
\end{equation}
in the limit $\lambda\to\infty$, where 
\begin{equation}\label{c12}\begin{split} 
C_1=&\frac{\Vol(X,g_p)\Vol(S^{n-1})}{n(2\pi)^n},\\ 
C_2=&\frac{\Vol(M,{a_0}^{1/2}h_0)\Vol(S^{n-1})}{2(2\pi)^n}. 
\end{split}\end{equation} 
If moreover we assume that $g_0$ is an exact cusp metric, then 
\[C_3=\frac{\Gamma\left(\tfrac{1-p}{2p}\right) } 
{2\sqrt{\pi}\Gamma\left(\tfrac{1}{2p}\right)} 
\zeta\left(\Delta_A^{h_0}, \frac{1}{p}-1\right),\] 
where $\Delta_A^{h_0}$ is the magnetic Laplacian on $M$ with 
potential $A_{|M}$ with
respect to the metric $h_0$ on $M$ defined in Section \ref{lcm}. 
\end{theorem}
We stress that the constants $C_1$ and $C_2$ do not
depend on the choice of $A$ or $B$, but only on the
metric. This fact provides some very interesting coupling constant effect, see Section
\ref{s:coupling}.

Note also that the hypotheses of the theorem are independent of
the choice of the vector potential $A$ inside the class of cusp $1$-forms.
Indeed, assume that $A'=A+dw$ for some $w\in\cun(X)$ is again a cusp 
$1$-form. Then $dw$ must be itself a cusp form, so
\[dw=x^2\px w \frac{dx}{x^2}+d^M w\in\cun(\oX,\ctsX).\]
Write this as $dw=\varphi'\frac{dx}{x^2}+\theta'_x$.
For each $x>0$ the form $\theta'_x$ is exact. By the Hodge decomposition 
theorem, the space of exact forms on $M$ is closed,
so the limit $\theta'_0$ is also exact. Now $dw$ is an exact cusp form,
in particular it is closed. This implies that $x^2\px\theta'=d^M \varphi'$.
Setting $x=0$ we deduce $d^M(\varphi')_{x=0}=0$, or equivalently
$\varphi'_{|x=0}$ is constant. Hence the conditions from Theorem \ref{t:thmag}
on the vector potential are satisfied simultaneously by $A$ and
$A'$.

\begin{proof}
From Theorem \ref{p:pmag}, the operator $x^{2p}\Delta_A$ is fully elliptic 
when $A$ is trapping. The result follows directly from 
\cite[Theorem 17]{wlom}. Let us explain the idea: the complex powers
$\Delta_A^{-s}$ belong to the cusp calculus, and are of trace class
for sufficiently large real part of $s$. The trace of the complex 
powers is holomorphic for such $s$ and extends meromorphically to $\cz$
with two families of simple poles, coming from the principal symbol and from 
the boundary. The leading pole governs eigenvalue asymptotics, by the Delange
theorem. In case the leading pole is double (by the superposition of the two 
families of poles), we get the logarithmic growth law.

The explicit computation of the constants (given by \cite[Prop.\ 14 and 
Lemma 16]{wlom}) is straightforward.
\end{proof}

\section{Analysis of the free case for non-trapping potentials}\label{s:free}

We have shown in Theorems \ref{p:pmag} and \ref{t:thmag} that 
the essential spectrum of the magnetic Laplacian $\Delta_A$ 
is empty when $A$ is trapping (see Definition \ref{def7}), 
and we have computed the asymptotics 
of the eigenvalues. We now consider the case of a non-trapping vector
potential $A$. One can guess that in this
last case, the essential spectrum is not empty 
when the metric is complete.
 In this section we concentrate on the unperturbed metric \eqref{gp'}
with the model non-trapping potential \eqref{acon}. We take
advantage of the decomposition in low- and high-energy functions
from Section
\ref{s:decomp}. The computation of the essential spectrum is based on 
Proposition \ref{p:free'} and on the diagonalization of the
magnetic Laplacian performed in Section \ref{s:diago}. In Section
\ref{s:conjop}, we construct a local conjugate operator and state the
Mourre estimate (Theorem \ref{t:mourre_0}). 

In Section
\ref{s:srlr} we introduce the classes of perturbation under which we
later give a limiting absorption principle. The perturbation
theory is developed in Section \ref{s:ntp}. We refer to Proposition
\ref{p:thema} for the question of the essential spectrum and to
Theorem \ref{t:mourre0} for its refined analysis under short/long
range perturbations. 

We localize the computation on the end
$X':=(0,\varepsilon) \times M\subset X$. We assume that
on $X':=(0,\varepsilon)\times M\subset X$ we have 
\begin{align}
g_p=&x^{2p}\left(\frac{dx^2}{x^4}+h_0\right),\label{gp'}\\
A_{\rm f}=&Cdx/x^2+\theta_0,\label{acon}
\end{align} 
where $C$ is a constant, $\theta_0$ is closed and independent of $x$, 
and the cohomology class $[\theta_0]\in H^1(M)$ is an integer multiple of
$2\pi$. By a change of gauge, one may assume that $C=0$. Indeed, it is enough
to subtract from $A_{\rm f}$ the exact $1$-form $d(-C/x)$.
\subsection{The high and low energy functions decomposition}\label{s:decomp}
Set $\dtt:=d^M+i\theta_0\wedge$. We now decompose the $L^2$ space as follows:
\begin{equation}\label{e:decomp}
L^2(X')=\Hr_{\rm l} \oplus \Hr_{\rm h},
\end{equation}
where $\Hr_{\rm l}:=\Kr \otimes \ker(\dtt)$ with 
$\Kr :=L^2\left((0,\varepsilon), x^{np-2} dx\right)$, and where
$\Hr_{\rm h}=\Kr\hat{\otimes} \ker(\dtt)^\perp$. We do not emphasize the
dependence on $\varepsilon$ for these spaces as the properties we are studying
are independent of $\varepsilon$. The subscripts ${\rm l},{\rm h}$ stand 
for low and high
energy, respectively. The importance of the \emph{low energy functions} space
$\Hr_{\rm l}$ is underlined by the following: 
\begin{proposition}\label{p:free'}
The magnetic Laplacian $\Delta_{A_{\rm f}}$ stabilizes the decomposition
\eqref{e:decomp} of $L^2(X')$.
Let $\Delta_{A_{\rm f}}^l$ and $\Delta_{A_{\rm f}}^h$ be the
Friedrichs extensions of the restrictions of $\Delta_{A_{\rm f}}$ to 
smooth compactly supported functions in $\Hr_{\rm l}$, respectively in
$\Hr_{\rm h}$. Then $\Delta_{A_{\rm f}}^h$ has compact resolvent, while
\begin{equation*}
\Delta^l_{A_{\rm f}}=\big(D^*D+c_0^2x^{2-2p}\big)\otimes 1,
\end{equation*} 
where $c_0:=((2-n)p-1)/2$, and $D:=x^{2-p}\px-c_0x^{1-p}$ acts in $\Krke$.
\end{proposition}
By Proposition \ref{p:cut}, the essential spectrum is
pres\-cri\-bed by the space of low energy functions. This kind of
decomposition can be found in various places in the literature,
see for instance \cite{lott} for applications to finite-volume
negatively-curved manifolds.

\begin{remark}\label{r:free'}
The dimension of the kernel of $\dtt$ equals the number of connected 
components of $M$ 
on which $A_f$ non-trapping. Indeed, take such a connected component $M_1$
and let $v$ be a primitive of $\theta_0$ on its universal cover. 
On different sheets of the cover, $v$ changes by $2\pi\zz$ so $e^{iv}$ 
is a well-defined function on $M_1$ which spans
$\ker(\dtt)$. This decomposition is also valid in the trapping
case, only that the low energy functions space is then $0$. Using
Proposition \ref{p:cut}, we obtain the emptiness of the essential spectrum for
the unperturbed metric and along the way a special case of Theorem
\ref{p:pmag} 
(we do not recover the full result in this way,
since the metric \eqref{pme} can \emph{not} be reached perturbatively
from the metric \eqref{gp'}). 
\end{remark} 

\proof 
We decompose the space of $1$-forms as the direct sum \eqref{decom}. 
Recall that $\delta_M$
is the adjoint of $d^M$ with respect to $h_0$. We compute
\begin{equation}\label{e:diam}\begin{split}
d_{A_{\rm f}}=&
\begin{bmatrix}
d^M+i\theta_0\wedge\\
x^2\px
\end{bmatrix}\\
d_{A_{\rm f}}^*=& x^{-np}\begin{bmatrix}\delta_M-i\theta_0\lrcorner& -x^2\px
\end{bmatrix}x^{(n-2)p}\\
\Delta_{A_{\rm f}}=& x^{-2p}(\dtt^*\dtt-(x^2\px)^2-(n-2)px (x^2\px)).
\end{split}\end{equation}
On the Riemannian manifold $(M,h)$, $\dtt^*\dtt$ is non-negative with discrete
spectrum. Since $-(x^2\px)^2-(n-2)px (x^2\px)=(x^2\px)^*(x^2\px)$ is non-negative, one has that
$\Delta_{A_{\rm f}}^h\geq \varepsilon^{-2p} \lambda_1$, where $\lambda_1$ is the
first non-zero eigenvalue of $\dtt^*\dtt$. By Proposition \ref{p:cut},
the essential spectrum is independent of
$\varepsilon$. By letting $\varepsilon \rightarrow 0$ we see that
it is empty; thus $\Delta_{A_{\rm f}}^h$ has compact resolvent. 
The assertion on $\Delta_{A_{\rm f}}^l$ is a straightforward computation.\qed 

\subsection{Diagonalization of the free magnetic Laplacian}\label{s:diago}
In order to analyze the spectral properties of $\Delta_{A_{\rm f}}$ on $(X,g_p)$, where
${A_{\rm f}}$ is given by \eqref{acon} with $C=0$ and $g_p$ by \eqref{gp'}, we go into some
``Euclidean variables''. We concentrate on the complete case, i.e.\ $p\leq 1$.
We start with \eqref{e:decomp} and work on $\Kr$. The first
unitary transformation is
\begin{align*}
L^2\left( x^{np-2} dx\right)\to L^2\left( x^{p-2} dx\right)&&
\phi\mapsto x^{(n-1)p/2}\phi.
\end{align*}
Then we proceed with the change of variables $z:=L(x)$, where $L$ is given by
\eqref{e:L}. Therefore, $\Kr$ is unitarily sent into $L^2\left((c,\infty),
dz\right)$ for a certain $c$. We indicate operators and spaces obtained in the
new variable with a subscript $0$.

Thanks to this transformation, we are pursuing our analysis on the
manifold $X_0=X$ endowed with the Riemannian metric 
\begin{align}\label{e:gp''}
dr^2+h,&&r\to\infty,\end{align}
on the end $X_0'=[1/2,\infty)\times M$. The subscript $r$ stands for radial.
The magnetic Laplacian is unitarily sent into an elliptic operator of
order $2$ denoted by $\Delta_0$. On $\cC^\infty_c(X_0')$, it acts by 
\begin{equation}\label{e:H_0}
\Delta_0:= Q_p\otimes {d_{\tilde\theta}}^*d_{\tilde\theta}+(-\partial^2_r +
V_p)\otimes 1, 
\end{equation}
on the completed tensor product $L^2([1/2,\infty), dr)\hat \otimes L^2(M, h)$,
where
\begin{align*}
V_p(r)=\begin{cases} ((n-1)/2)^2\\c_0r^{-2}\end{cases} 
&& \mbox{and } Q_p(r)=\begin{cases}
e^{2r}\\
((1-p)r)^{2p/(1-p)}
\end{cases} &&\mbox{for } \begin{cases}
 p=1\\
 p<1.
\end{cases}
\end{align*}
Recall that $c_0$ is defined in Proposition \ref{p:free'} and
$d_{\tilde\theta}=d_M+i\theta_0\wedge$. 

We denote also by $L_0$ the operator of multiplication corresponding
to $L$, given by \eqref{e:L}, in the new variable $r$. It is bounded from
below by a positive constant, equals $1$ on the compact part and on
the trapping ends, and equals $r\mapsto r$ on the non-trapping ends.

Let $\Hr^s_0:=\Dc((1+\Delta_0)^{s/2})$ for $s>0$. By identifying $\Hr_0$ with
$\Hr^*_0$ by the Riesz isomorphism, by duality, we define $\Hr^s_0$ for $s<0$
with $\Hr^{-s}_0$. We need the next well-known fact.

\begin{lemma}\label{l:cchi}
For every $\gamma\in\cC^\infty_c(X)$, we have that 
$\gamma:\Hr^s\subset\Hr^s$ for all $s\in\R$
\end{lemma}
\proof
A computation gives that there is $c$ such that
$\|(\Delta_0+i)\gamma\varphi\|\leq c \|(\Delta_0+i)\varphi\|$, for all
$\varphi\in\cC^\infty_c(X_0)$. Since $\cC^\infty_c(X_0)$ is a core for $\Delta_0$,
we get the result for $s=2$. By induction we get it for $s\in 2\N$. Duality
and interpolation give it for $s\in\R$.\qed

\subsection{The local conjugate operator}\label{s:conjop}
In this section, we construct some conjugate operators in order to
establish a Mourre estimate. By mimicking the case of the Laplacian, see
\cite{FH}, one may use the following localization of the generator of
dilations. Let $\xi\in\cC^\infty([1/2,\infty))$ such that the support of
$\xi$ is contained in $[2,\infty)$ and that $\xi(r)=r$ for $r\geq 3$
and let $\tilde \cchi\in\cC^\infty([1/2,\infty))$ with
support in $[1,\infty)$, which equals $1$ on $[2,\infty)$. By abuse of
notation, we denote $\tilde\cchi\otimes 1\in\cC^\infty(X_0)$ with
the same symbol. Let $\cchi:=1-\tilde\cchi$. On $\cC^\infty_c(X_0)$
we set: 
\begin{equation}\label{e:S0}
S_\infty:= \left( -i (\xi \partial_r +\partial_r \xi)\otimes P_0\right) \tilde\cchi,
\end{equation} 
where $P_0$ is the orthogonal projection onto $\ker(\dtt)$. 
The presence of $P_0$ comes from the decomposition in low and high energies.
\begin{remark}\label{r:sa}
By considering a $C_0$-group associated to a vector field on $\R$ like in
\cite[Section 4.2]{ABG}, one shows that $-i(\xi \partial_r +\partial_r \xi)$
is essentially self-adjoint on $\cC^\infty_c(\R)$. This $C_0$-group acts
trivially away from the support of $\tilde{\cchi}$ and then it is easy to
construct another $C_0$-group $G_0$ which acts like the first on
$L^2(\supp(\tilde{\cchi}))\otimes\ker(\dtt)$, and trivially 
on the rest of $L^2(X_0)$. Let $\overline{S_\infty}$ be 
the generator of $G_0$. Since $G_0(t)$
leaves invariant $\cC^\infty_c(X_0)$ the Nelson lemma implies that
$\overline{S_\infty}$ is the closure of $S_\infty$, in other words
$\overline{S_\infty}$ is essentially self-adjoint on $\cC^\infty_c(X_0)$. This
kind of approach has been used in \cite{bouclet} for instance. 
\end{remark}

Therefore, we denote below also by $S_\infty$
the self-adjoint closure of the operator
defined by \eqref{e:S0} on $\cC^\infty_c(X_0)$.

However this operator is \emph{not} suitable for very singular perturbations
like that of the metric considered in this paper. To solve this problem,
one should consider a conjugate operator more ``local in energy''.
Concerning the Mourre estimate, as it is local in energy for the Laplacian,
one needs only a conjugate operator which fits well on this level of energy.
Considering singular perturbation theory, the presence of differentials in
\eqref{e:S0} is a serious obstruction; the idea is to replace the conjugate
operator with a multiplication operator in the analysis of perturbations,
therefore reducing the r\^ole of derivatives within it. The approach has
been used for Dirac operators for instance in \cite{GM} to treat very singular
perturbations. The case of Schr\"odinger operators is
summarized in \cite[Theorem 7.6.8]{ABG}. We set:
\begin{equation}\label{e:SR}
S_R:=\tilde\cchi \left(\big(\Phi_R(-i\partial_r) \xi +\xi
\Phi_R(-i\partial_r)\big)\otimes P_0\right) \tilde\cchi,
\end{equation}
where $\Phi_R(x):=\Phi(x/R)$ for some $\Phi\in\cC^\infty_c(\R)$ satisfying 
$\Phi(x)=x$ for all $x\in[-1,1]$. The operator $\Phi_R(-i\partial_r)$ is 
defined on $\rz$
by $\Fr^{-1}\Phi_R\Fr$, where $\Fr$ is the unitary Fourier transform. Let us
also denote by $S_R$ the closure of this operator.

Unlike \eqref{e:S0}, $S_R$ does not stabilize $\cC^\infty_c(X_0)$
because $\Phi_R(-i\partial_r)$ acts like a convolution with a function with
non-compact support. This subspace is sent
into $\tilde\cchi \Sr(\R)$, where $\Sr(\R)$ denotes the Schwartz space. 
To motivate the subscript $R$, note that $S_R$ tends strongly in the
resolvent sense to $S$, as $R$ goes to infinity. We give some
properties of $S_R$. The point (\ref{p2}) is essential to be able to replace
$S_R$ by $L_0$ in the theory of perturbations.
The point (\ref{p3}) is convenient to be able to express a limiting absorption
principle in terms of $L_0$, which is very explicit. Of course
these two points are false for $S_\infty$ and this explains
why we can go further in the perturbation theory
compared to the standard approach.

\begin{lemma}\label{l:sa}
Let $S_R$ denote the closure of the unbounded operator \eqref{e:SR}.
\begin{enumerate}
\item For all $R\in[1,\infty]$, the operator $S_R$ is essentially
self-adjoint on $\cC^\infty_c(X)$.
\item For $R$ finite, $L_0^{-2}S_R^2: \cC^\infty_c(X_0)
\rightarrow \Dc(\Delta_0)$ extends to a bounded operator in $\Dc(\Delta_0)$.
\label{p2}
\item For $R$ finite, $\Dc(L_0^s)\subset \Dc(|S_R|^s)$ for all $s\in[0,2]$.
\label{p3}
\end{enumerate}
\end{lemma}
\proof The case $R=\infty$ is discussed in Remark \ref{r:sa}, so assume that 
$R$ is finite. We compare $S_R$ with $L_0$, defined in Section \ref{s:diago},
which is essentially self-adjoint on $\cC^\infty_c(X)$. Noting that it
stabilizes the decomposition \eqref{e:decomp}, we write also by
$L_{0,\rm l}$ its
restriction to $\Hr_{\rm l}$, which is simply the multiplication by
$r$. On $\cC^\infty_c(X_0)$,
\begin{equation*}
S_R= \tilde\cchi\big(2\Phi_R(-i\partial_r)\xi L_{0,\rm l}^{-1}+
[\xi, \Phi_R(-i\partial_r)] L_{0,\rm l}^{-1} \big)\otimes P_0 \tilde\cchi L_0
\end{equation*}
Noting that $\xi L_{0,\rm l}^{-1}$ is bounded and that $\xi'\in
L^\infty$ and using Lemma \ref{l:pseudo}, we get $\|S_R \varphi\|\leq a
\|L_0 \varphi\|$, for all $\varphi\in\cC^\infty_c(X_0)$.

On the other hand, $[S_R, L_0]$
is equal to the bounded operator
\begin{align*}
[S_R, L_0]=& \tilde\cchi\big([\Phi_R(-i\partial_r), L_{0,\rm l}]\xi
L_{0,\rm l}^{-1}\big) \otimes P_0 \tilde\cchi L_0\\ &
+ L_0 \tilde\cchi\big( L_{0,\rm l}^{-1} \xi [\Phi_R(-i\partial_r), L_{0,\rm l}] \big) \otimes P_0\tilde\cchi. 
\end{align*}
This gives $|\langle S_R \varphi, L_0\varphi \rangle -\langle L_0\varphi,
S_R\varphi \rangle| \leq b \|L_0^{1/2}\varphi\|^2$, for all
$\varphi\in\cC^\infty_c(X_0)$. Finally, one uses \cite[Theorem X.37]{RS2} to
conclude that $S_R$ is essentially self-adjoint.

On $\cC^\infty_c(X_0)$, we have 
\begin{equation}\label{e:LS}
 L^{-2}_0S_R^2= (2\cchi L^{-1}_{0,\rm l}\xi
\Phi_R(-i\partial_r)\cchi\otimes P_0 +
\cchi L^{-1}_{0,\rm l} [\Phi_R(-i\partial_r), \xi]\cchi\otimes P_0)^2.
\end{equation}
All these terms are bounded in $L^2(X_0)$ by Lemma \ref{l:pseudo} and
by density. We now compute $\Delta_0 L^{-2}_0S_R^2$. By Lemma \ref{l:cchi},
it is enough to show that $\Phi_R(-i\partial_r)$ and
$[\Phi_R(-i\partial_r), \xi]$ stabilize the domain $\Delta$ in
$L^2(\R)$. The first one commutes with $\Delta$. For the second one,
we compute on $\cC^\infty_c(\R)$. Since $[\Phi_R(-i\partial_r), \xi]$
is bounded in $L^2(\R)$, it is enough to show that
the commutator $[\Delta,[\Phi_R(-i\partial_r), \xi]]$ is also bounded
in $L^2(\R)$. By Jacobi's identity, it is equal to
$[\Phi_R(-i\partial_r), [\Delta, \xi]] = [\Phi_R(-i\partial_r), 2
\xi'\partial_r + \xi'']= 2\Phi_R(-i\partial_r)\partial_r \xi'-2
\Phi_R(-i\partial_r) \xi'' + 2 \xi'\Phi_R(-i\partial_r) +\xi''$. This
is a bounded operator in $L^2(\R)$ and we get point (2).

We now note that \eqref{e:LS} is bounded in $L^2(X_0)$. Then, since
$S_R^2L_0^{-2}$ is also bounded, we get $\|S_R^2\varphi\|^2 \leq c
\|L^2_0\varphi\|$ for all $\varphi \in \cC^\infty(X_0)$. Taking a
Cauchy sequence, we deduce $\Dc(L^2_0)\subset \Dc(S_R^2)$. An
argument of interpolation gives point (3).\qed

The aim of this section is the following Mourre estimate. 

\begin{theorem}\label{t:mourre_0}
Let $R\in [1,\infty]$. Then $e^{itS_R} \Hr^2_0\subset \Hr^2_0$ and
$\Delta_0\in \cC^{2}(S_R, \Hr^2_0, \Hr_0)$. Given an
interval $\cJ$ inside $\sigma_{\rm ess}(\Delta_0)$, there exist
$\varepsilon_R> 0$ and a compact operator $K_R$ such that 
\begin{equation*}
E_\cJ(\Delta_0) [\Delta_0,iS_R] E_\cJ(\Delta_0) \geq (4\inf(\cJ)
- \varepsilon_R) E_\cJ(\Delta_0)+ K_R
\end{equation*} 
holds in the sense of forms, and such that $\varepsilon_R$ tends to $0$
as $R$ goes to infinity.
\end{theorem} 
\proof The regularity assumptions follow from Lemmata \ref{l:regu0} and
\ref{l:regu}.
The left hand side of \eqref{e:mourre_0} is the commutator $[\Delta_0,
iS_R]$ in the sense of forms. It extends to a bounded operator in 
$\Bc(\Hr^2_0, \Hr^{-2}_0)$ since $\cC^\infty_c(X_0)$ is a core for
$\Delta_0$. We can then apply the spectral measure and obtain the
inequality using Lemma \ref{l:mourre_0}.\qed 

Compared to the method from \cite[Lemma 2.3]{FH} (for the case
$R=\infty$), we have a relatively more direct proof based on Lemma
\ref{l:mourre_0}. However this has no real impact on applications of
the theory.

We now go in a series of lemmata to prove this theorem. Given a
commutator $[A,B]$, we denote its closure by $[A,B]_0$. 
\begin{lemma}\label{l:commu}
For $R=\infty$, the commutators 
$[\Delta_0,iS_\infty]_0$ and $[[\Delta_0,iS_\infty],iS_\infty]_0$
belong to $\Bc(\Hr^2_0, \Hr_0)$. For $R$ finite, $[\Delta_0,iS_R]_0$ and
$[[\Delta_0,iS_\infty],iS_R]_0$ belong to $\Bc(\Hr_0)$. Moreover, if
$p=1$, all higher commutators extend to bounded operators in $\Bc(\Hr^2_0,
\Hr_0)$ for $R=\infty$ and in $\Bc(\Hr_0)$ for $R<\infty$.
\end{lemma}

\proof Let $\varphi\in\cC^{\infty}_c(X_0')$ such that
$\varphi=\varphi_\R\otimes \varphi_M$ where
$\varphi_M\in\cC^\infty(M)$ and $\varphi_\R\in\cC^\infty_c([1/2,
\infty))$. Note that $P_0\varphi_M\in\cC^\infty(M)$ by the Hodge decomposition,
since $d_A^2=0$. Applying the brackets to $\varphi$, by a straightforward
computation, we get
\begin{equation}\label{e:commu}
[\Delta_0,iS_\infty]=-(4\xi'\partial^2_r+4\xi''\partial_r+\xi'''-2V_p')\tilde
\cchi\otimes P_0.
\end{equation}
By linearity and density, we get $\|[\Delta_0,iS_\infty]\varphi\|\leq
C\|(\Delta_0+i)\varphi\|$ for all $\varphi\in\cC^\infty_c(X_0')$. Take now
$\varphi\in\cC^\infty_c(X_0)$, considering the support of the
commutator, we get
\begin{eqnarray*}
\|[\Delta_0,iS_\infty]\varphi\|&=& \|[\Delta_0,iS_\infty]\widetilde\Xi\varphi\|
\\ 
&\leq& C\|(\Delta_0+i) \varphi\| + \|[\Delta_0, \widetilde\Xi] \varphi\|
\leq C'\|(\Delta_0+i) \varphi\|.
\end{eqnarray*}
where $\widetilde\Xi\in\cC^{\infty}(X)$ with support in $X_0'$ such that
$\widetilde\Xi|_{[1,\infty)\times M}=1$.
Therefore, since $\cC^\infty_c(X_0)$ is a core for $\Delta_0$,
we conclude $[\Delta_0,iS_\infty]_0\in\Bc(\Hr^2_0, \Hr_0)$.

In the same way we compute
\begin{eqnarray}\nonumber
[[\Delta_0,iS_\infty], iS_\infty]&=& -\big((16(\xi')^2-8\xi\xi'')\partial^2_r + (24
\xi'\xi''-8\xi\xi''')\partial_r \\\label{e:commu2}
&&+4 (\xi'')^2+4\xi'''-2 \xi'''' - 4V_p''\big)\tilde\cchi\otimes P_0.
\end{eqnarray}
and get $[[\Delta_0,iS_\infty],iS_\infty]_0\in\Bc(\Hr^2_0, \Hr_0)$. For $p=1$, the boundedness of higher commutators follows easily by induction ($V_p=0$ in this case).

We compute next the commutators of $\Delta_0$ with $S_R$. 
As above we compute for $\varphi=\varphi_r\otimes\varphi_M$. For brevity, 
we write $\Phi_R$ instead of $\Phi_R(-i\partial_r)$.

As $\Phi_R$ is not a local operator, we first note that the
commutator $[\Delta_0, S_R]$ 
could be taken in the operator sense. Indeed, $\tilde\cchi$ sends
$\varphi_r$ to $\cC^\infty_c(\R)$ (note that $[1/2,\infty)$ is
injected in a canonical way into $\R$), then $\Phi_R \xi +\xi\Phi_R$
sends to the Schwartz space $\Sr(\R)$ and finally $\tilde\cchi$ sends
to $\tilde\cchi\Sr(\R)$ which belongs to $\Dc(\Delta_0)$.

We compute $[\partial_r^2, \tilde\cchi \big(\Phi_R \xi +\xi
\Phi_R\big) \tilde\cchi]\otimes P_0$. Against $\varphi_r\otimes \varphi_M$, 
we have:
\begin{eqnarray}\nonumber
[\partial_r^2, \tilde\cchi \Phi_R \xi +\xi
\Phi_R \tilde\cchi]&=& \Xi [\partial_r^2, \tilde\cchi \Phi_R \xi +\xi
\Phi_R \tilde\cchi] \Xi = \Xi [\partial_r^2, \Phi_R r + r\Phi_R ] \Xi
+\Psi_\comp \\ \label{e:commu5}
&=& 4\Xi \partial_r\Phi_R \Xi +\Psi_\comp = 4\tilde\cchi \partial_r\Phi_R
\tilde\cchi +\Psi_\comp,
\end{eqnarray}
where $\Psi_\comp$ denotes a pseudo-differential operator with compact
support such that its support in position is in the interior of
$X_0'$. For $p<1$, the potential part $V_p$ arises. We 
treat its first commutator: 
\begin{eqnarray} \nonumber
[V_p, \tilde\cchi \Phi_R \xi +\xi \Phi_R\tilde\cchi]
&=&\tilde\Xi[V_p, 2\Phi_R r -i\Phi'_R]\tilde\Xi +\Psi_\comp
\\ \label{e:commu6}
&=&\tilde\cchi\big(2[V_p,\Phi_R]\xi-i[V_p, \Phi_R'] \big)\tilde\cchi+\Psi_\comp,
\end{eqnarray}
where $\Phi_R'=\Phi_R'(-i\partial_r)$. Applying Lemma \ref{l:pseudo}, we get that $[V_p, \Phi_R]\xi$ and $[V_p, \Phi_R']$ are bounded in $L^2(\R)$ also. Therefore, using like above $\tilde{\Xi}$, we get $\|[\Delta_0,iS_R]\varphi\|\leq C\|\varphi\|$ for all $\varphi\in\cC^\infty_c(X_0)$. This implies that $[\Delta_0,iS_R]_0\in\Bc(\Hr_0)$.

For higher commutators, the $n$-th commutator with $\tilde\cchi
\Phi_R \xi +\xi \Phi_R \tilde\cchi$ is given by $2^n\tilde\cchi
\partial_r^n\Phi_R \tilde\cchi +\Psi_\comp$. Note that $\partial_r^n\Phi_R$ is a compactly supported function of $\partial_r$, so the contribution of this term is always bounded.

Consider now the second commutator of $V_p$. As above, since we work up to 
$\tilde\Xi$, $\tilde{\cchi}$ and $\Psi_\comp$, it is enough to show that the next commutator defined on $\Sr(\R)$ extend to bounded operators in $L^2(\R)$. We treat only the most singular part of the second commutator:
\begin{eqnarray}\nonumber
[[V_p, \Phi_R]r, \Phi'_R r]&=& [[V_p, \Phi_R], \Phi_R r]r + [V_p, \Phi_R][r,\Phi_R]r
\\ \nonumber
&\hspace*{-4cm }=& \hspace*{-2cm } [[V_p, \Phi_R], \Phi_R ]r^2 + \Phi_R [[V_p, \Phi_R], r]r -i [V_p, \Phi_R]\Phi_R'r
\\ \nonumber
&\hspace*{-4cm }=& \hspace*{-2cm } [[V_p, \Phi_R], \Phi_R ]r^2 - \Phi_R [[r, \Phi_R], V_p]r -i [V_p, \Phi_R][\Phi_R',r] -i [V_p, \Phi_R]r \Phi_R'
\\ \nonumber
&\hspace*{-4cm }=& \hspace*{-2cm } [[V_p, \Phi_R], \Phi_R ]r^2 + i \Phi_R [\Phi_R', V_p]r + [V_p, \Phi_R] \Phi_R'' -i [V_p, \Phi_R]r \Phi_R'
\end{eqnarray}
These terms extend to bounded operators by Lemma \ref{l:pseudo}. \qed

The following lemma is the key-stone for the Mourre estimate.
\begin{lemma}\label{l:mourre_0}
For all $R\in[1,\infty]$, there exists $K\in\Kc(\Hr^2_0, \Hr^{-2}_0)$
and $N_R\in\Bc(\Hr^2_0, \Hr^{-2}_0)$ such that
\begin{eqnarray}\nonumber
\langle \Delta_0\varphi, iS_R\varphi \rangle+ \langle iS_R
 \varphi, \Delta_0\varphi \rangle &=& 
\\ \label{e:mourre_0}
&& \hspace*{-4cm}4\langle \varphi,
 (\Delta_0-\inf(\sigma_{\rm ess}(\Delta_0)) \varphi\rangle + \langle \varphi, N_R
 \varphi\rangle +\langle \varphi, K\varphi\rangle, 
\end{eqnarray}
for all $\varphi\in\cC^\infty_c(X_0)$ and such that $\|N_R\|_{\Bc(\Hr^2_0,
\Hr^{-2}_0)}$ tends to $0$ as $R$ goes to infinity.
\end{lemma}

\proof
First note that the essential spectrum of $\Delta_0$ is $[V_p(\infty), \infty)$, by
Proposition \ref{p:thema}. We now act in three steps.
Let $\Xi$ be like in the proof of Lemma \ref{l:commu} and let
$\varphi\in\cC^\infty_c(X_0)$. Since $\tilde\cchi$ and $\Xi$ have disjoint 
supports, one has for all $R$ that $\langle \varphi,
[\Delta_0,iS_R] \varphi\rangle= \langle\widetilde\Xi \varphi, [\Delta_0,iS_R]
\widetilde\Xi\varphi\rangle$.

For the first step, we start with $R=\infty$. 
By \eqref{e:commu} and since $\xi'=1$ on
$[2,\infty)$, the Rellich-Kondrakov lemma gives
\begin{equation}\label{e:Mstep1}
\langle \widetilde\Xi \varphi, [\Delta_0,iS_\infty] \widetilde\Xi\varphi\rangle
=\langle \widetilde\Xi \varphi, 4(\Delta_0 -V_p(\infty))(1\otimes
P_0)\widetilde\Xi \varphi \rangle+ \langle \varphi, K_1 \varphi \rangle 
\end{equation}
for a certain $K_1\in\Kc(\Hr^2_0,\Hr_0^{-2})$. Indeed,
$1-\xi', V_p-V_p(\infty), \xi'',\xi''', V_p'$ belong to
$\Kc(\Hr^1_0,\Hr_0)$ since they tend to $0$ at infinity.

We consider now $R$ finite. We add \eqref{e:commu5} and
\eqref{e:commu6}. We have 
\begin{equation}\label{e:Mstep1'}
\langle \widetilde\Xi \varphi, [\Delta_0,iS_R] \widetilde\Xi\varphi\rangle
=\langle \widetilde\Xi \varphi, 4(\Delta_0 -V_p(\infty) - T_R)(1\otimes
P_0)\widetilde\Xi \varphi \rangle
+ \langle \varphi, K_2 \varphi \rangle 
\end{equation}
for a certain $K_2=K_2(R)\in\Kc(\Hr^2_0,\Hr_0^{-2})$ and with $T_R=
\partial_r(\partial_r-\Phi_R(\partial_r))$. The compactness of $K_2$ follows by
noticing that $L_0^{-1}\in \Kc(\Hr^2_0,\Hr_0)$ and that 
$L_0[V_p, i S_R]_0\in \Bc(\Hr_0,\Hr_0)$, by Lemma \ref{l:pseudo}. 
We control the size of $T_R$ by showing that
$\|\widetilde\Xi T_R (1\otimes P_0)\widetilde\Xi \|_{\Bc(\Hr^2_0,
\Hr^{-2}_0)}$ tends to $0$ as $R$ goes to infinity. By Lemma
\ref{l:cchi}, one has that $\widetilde\Xi$ stabilizes $\Hr^{\pm
2}_0$, therefore $-\partial_r^2 \widetilde\Xi$ belongs to
$\Bc(\Hr^2_0, L^2(\R))$. It remains to note that $(-\partial_r^2+i)^{-2}T_R$
tends to $0$ in norm by functional calculus, as $R$ goes to infinity.

The second step is to control the high energy functions part. Consider the
Friedrichs extension of $\Delta_0$ on $\Hr_{0,h}:=L^2([1/2,
\infty))\otimes P_0^\perp L^2(M)$. We have:
\begin{equation}\label{e:Mstep2}
\langle \widetilde\Xi \varphi, \Delta_0 P_0^\perp\widetilde\Xi\varphi\rangle=
\langle (\Delta_0 P_0^\perp+i)^{-1}(\Delta_0 P_0^\perp+i)
\widetilde\Xi \varphi, \Delta_0 P_0^\perp\widetilde\Xi\varphi\rangle
=\langle \varphi, K_1 \varphi\rangle 
\end{equation}
where $K_1\in\Kc(\Hr_0^2, \Hr_0^{-2})$.
Indeed, note first that $\Delta_0 P_0^\perp\tilde\cchi \in
\Bc(\Hr_0^2,\Hr_{0,h})$ and that $(\Delta_0 P_0^\perp+i)^{-1}
\in\Kc(\Hr_{0,h})$, since $Q_p$ in \eqref{e:H_0} 
goes to infinity. Therefore the left hand side belongs $\Kc(\Hr_0^2,\Hr_0)$
and the right hand side belongs to $\Bc(\Hr^2_0,\Hr_0)$. 

The third step is to come back on the whole manifold. 
It is enough to note that $[\Delta_0,\tilde\cchi]\in\Kc(\Hr^1_0,
\Hr_0)$ and to add \eqref{e:Mstep1}, \eqref{e:Mstep1'} with
\eqref{e:Mstep2}.\qed 

We now turn to the regularity assumptions. Lemma \ref{l:C1} plays a
central r\^ole. 

\begin{lemma}\label{l:regu0}
For $R\in [1, \infty]$, one has $\Delta_0\in\cC^1(S_R)$ and 
$e^{itS_R} \Dc(\Delta_0) \subset\Dc(\Delta_0)$. 
\end{lemma}
\proof
We start by showing that $\Delta_0\in\cC^1(S_R)$. We check the hypothesis of Lemma
\ref{l:C1}. Let $\cchi_n(r):=\cchi(r/n)$ and $\Dr=\cC^\infty_c(X_0)$. Remark
that $\supp(\cchi_n')\subset[n,2n]$ and that $\xi \cchi_n^{(k)}$ tends
strongly to $0$ on $L^2(\R^+)$, for any $k\geq 1$. By the uniform boundedness 
principle, this implies that $\sup_n\|\cchi_n\|_{\Dc(H)}$ is finite.

Remark \ref{r:sa} and Lemma \ref{l:sa} give that $\Dr$ is a core for
$S_R$. Assumption (1) is obvious, assumption (2) holds since $(1-\cchi_n)$ has
support in $[2n,\infty)$ and assumption (3) follows from the fact that
$H$ is elliptic, so the resolvent of $\Delta_0$ sends $\Dr$ into $\cC^\infty(X_0)$. 
The point \eqref{e:C1} follows from Lemma \ref{l:mourre_0}.
We now show that \eqref{e:comm} is true. Let
$\phi\in\cC^\infty(X_0)\cap\Dc(\Delta_0)$. 
We have $[\Delta_0,\cchi_n]\phi=[\Delta_0,\cchi_n]\tilde\cchi\phi=2
\cchi_n'\partial_r\tilde\cchi\phi + \cchi_n''\tilde\cchi\phi$. We
have $iS_\infty[\Delta_0,\cchi_n]\phi=
2\xi\cchi_n'\partial^2_r P_0\tilde \cchi\phi +(4\cchi\xi_n''
+\xi'\cchi_n')\partial_r P_0 \tilde\cchi\phi+
(2\xi\cchi_n'''+\xi'\cchi_n'')P_0 \tilde\cchi\phi$ and, for a finite $R$, we
get $iS_R[\Delta_0,\cchi_n]\phi=
\tilde\cchi(2\Phi_R(\partial_r)\xi + [\xi, \Phi_R(\partial_r)])
(2\cchi_n'\partial_r P_0\tilde\cchi\phi + \cchi_n'' P_0\tilde\cchi\phi)$.
Both terms are tending to $0$ because of the previous remark, Lemma
\ref{l:cchi} and the fact that $[\xi, \Phi_R(\partial_r)]$ is
bounded by Lemma \ref{l:pseudo}. From that, we can apply the lemma
and obtain $H\in\cC^1(S_R)$. 

 By Lemma \ref{l:commu}, we have that $[\Delta_0, iS_R]_0\in\Bc(\Hr_0^2,\Hr_0)$ and
\cite[Lemma 2]{GG0} gives that $e^{itA}\Hr_0^2\subset\Hr_0^2$. \qed 

The invariance of the domain under the group $e^{itS_R}$ implies that
$e^{itS_R}\Hr^s_0\subset\Hr^s_0$ for $s\in [-2, 2]$ by duality and
interpolation. This allows one to define the class $\cC^k(S_R,
\Hr^{s}_0,\Hr^{-s}_0)$ for $s\in [-2,2]$, for instance; we recall that a
self-adjoint operator $H$ is in this class if $t\mapsto e^{itS_R}He^{-itS_R}$
is strongly $C^k$ from $\Hr^{s}_0$ to $\Hr^{-s}_0$.

\begin{lemma}\label{l:regu}
Let $R\in [1, \infty]$. Then $\Delta_0$ belongs to $\cC^2(S_R, \Hr^2_0,\Hr_0)$ for $p\leq 1$.
\end{lemma} 
\proof
From Lemma \ref{l:commu}, the commutators with $S_R$ extend to bounded 
operators from $\Hr^2_0$ to $\Hr_0$. \qed

We finally give an estimation of commutator that we have used above.
\begin{lemma}\label{l:pseudo}
Let $f\in\cC^{0}(\R)$ with polynomial growth, $\Phi_j\in\cC^\infty_c(\R)$ and
$g\in\cC^k(\R)$ with bounded derivatives. Let $k\geq 1$. Assume that
$\sup_{t\in\R, |s-t|\leq 1}|f(t)g^{(l)}(s)|<\infty$, for all $1\leq l\leq k$. 
Then the operator 
$f[\Phi_1(-i\partial_r), [\Phi_2(-i\partial_r)\ldots 
[\Phi_k(-i\partial_r), g]\ldots]$, defined on $\cC^\infty_c(\R)$, 
extends also to a bounded operator.
\end{lemma}

\proof Take $k=1$. We denote with a hat the unital Fourier transform. We get
\begin{equation*}
 (f[\Phi(-i\partial_r), g]\varphi)(t)= \frac{1}{2\pi}\int
 \widehat{\Phi}(s-t)f(t)(g(s)-g(t))\varphi(s) ds,
\end{equation*}
for $\varphi\in\cC^\infty_c(\R)$. In order to show that the $L^2$ norm of the
left hand side is uniformly bounded by $\|f\|_2$, we separate the integral in
$|s-t|$ lower and bigger than one. We start with the first part. Recalling
$t^k\widehat\Phi(t)=(-1)^k\widehat{\Phi^{(k)}}(t)$, 
$\Phi$ is replaced by $\Phi'$ when we divide by $|s-t|$ the term with
$g$. Now since $\widehat{\Phi'}\in L^1$, $\sup_{t\in\R, |s-t|\leq
 1}|f(t) (g(s)-g(t))/ (s-t)|$ is finite, and the convolution by a $L^1$
function is bounded in $L^2$, we control this part of the integral. We turn to
the part $|s-t|\geq 1$. Let $R\in\N$ such that $|f(t)|\leq C(1+ |t|^R)$.
We let appear the Fourier transform of $\Phi^{(R+1)}$. Now to conclude,
note that $\sup_{s,t\in\R}|(g(s)-g(t))/ (s-t)|$ is finite and that
$t\mapsto|f(t)/(s-t)^R|$ is in $L^\infty$ uniformly in $s$. For higher $k$,
one repeats the same decomposition and let appear the terms in
the $l$-th derivative of $g$ by regrouping terms.\qed

\subsection{A short-range and long-range class of perturbations}\label{s:srlr}
In the early versions of Mourre theory, one asked
$[[H,S_R], S_R]$ to be $H-$bounded to obtain refined results of the
resolvent like the limiting absorption principle and the H\"older
regularity of the resolvent. In this section, we check the
optimal class of regularity $\cC^{1,1}(S_R)$, for $R$ finite. This is
a weak version of the two-commutators hypothesis. We refer to
\cite{ABG} for definition and properties. This is the optimal 
class of operators which give a limiting absorption principle for $H$ in 
some optimal Besov spaces associated to the conjugate operator $S_R$. 

The operator $\Delta_0$ belongs to $\cC^2(S_R)$ by Lemma \ref{l:regu} and
therefore also to $\cC^{1,1}(S_R)$. We now consider perturbations of
$\Delta_0$ which are also in $\cC^{1,1}(S_R)$. We define two classes.

Consider a symmetric differential operator $T:
\Dc(\Delta_0)\rightarrow \Dc(\Delta_0)^*$. Take $\theta_{\rm
sr}\in\cC_c^\infty((0,\infty))$ not identically $0$; $V$ is said to be
\emph{short-range} if 
\begin{equation}\label{e:sr}
\int_1^\infty \Big\|\theta_{\rm sr}\Big(\frac{L_0}{r}\Big) T \Big\|_{\Bc(\Dc(\Delta_0),
\Dc(\Delta_0)^*)} dr <\infty.
\end{equation}
and to be \emph{long-range} if
\begin{align}
\int_1^\infty & \Big\|[T, L_0]\theta_{\rm lr}\Big(\frac{L_0}{r}\Big)
\Big\|_{\Bc(\Dc(\Delta_0), \Dc(\Delta_0)^*)} 
+ \Big\|\tilde\Xi[T, P_0]L_0\theta_{\rm
lr}\Big(\frac{L_0}{r}\Big)\tilde\Xi\Big\|_{\Bc(\Dc(\Delta_0),
\Dc(\Delta_0)^*)}\nonumber
\\
&+ \Big\|\tilde\Xi[T, \partial_r]P_0 L_0\theta_{\rm
lr}\Big(\frac{L_0}{r}\Big)\tilde\Xi \Big\|_{\Bc(\Dc(\Delta_0),
\Dc(\Delta_0)^*)} \frac{dr}{r} <\infty.\label{e:lr}
\end{align}
where $\theta_{\rm lr}$ is the characteristic function of $[1,\infty)$
in $\R$ and where $\widetilde\Xi\in\cC^{\infty}(X)$ with support in
$X_0'$ such that $\widetilde\Xi|_{[1,\infty)\times M}=1$. 

The first condition is evidently satisfied if there is $\varepsilon
\pg 0$ such that 
\begin{equation}\label{e:sr'}
	\|L^{1+\varepsilon}_0 T\|_{\Bc(\Dc(\Delta_0),
\Dc(\Delta_0)^*)} \pp \infty
\end{equation}
and the second one if
\begin{align}\nonumber
\|L^{\varepsilon}_0 [T,L_0]\|_{\Bc(\Dc(\Delta_0), \Dc(\Delta_0)^*)}	&+
\|L^{1+\varepsilon}_0
\widetilde\Xi[T,P_0]\widetilde\Xi\|_{\Bc(\Dc(\Delta_0),
\Dc(\Delta_0)^*)} 
\\ \label{e:lr'}
&+ \|L^{1+\varepsilon}_0
\widetilde\Xi[T,\partial_r]P_0\widetilde\Xi\|_{\Bc(\Dc(\Delta_0),
\Dc(\Delta_0)^*)}\pp \infty 
\end{align}
The condition with $P_0$ essentially tells that the non-radial part of $T$ 
is a short-range perturbation. This is why we will ask the long-range
perturbation to be radial. To show that the first class is in 
$\cC^{1,1}(S_R)$ for $R$ finite, one use \cite[Theorem 7.5.8]{ABG}. 
The hypotheses are satisfied thanks to Lemmata \ref{l:sa} and
\ref{l:L}. Concerning the second class, one shows that $[T,
S_R]\in\cC^{0,1}(S_R)$ by using \cite[Proposition 7.5.7]{ABG} (see the
proof of \cite[Proposition 7.6.8]{ABG} for instance). 

We go back to the $x$ coordinate.
For $\Gr=\Dc(\Delta^{1/2})$, the short
and long-range perturbation of the electric and magnetic perturbations
are given by: 

\begin{lemma}\label{l:magne}
Let $V\in L^\infty(X)$ and $\tilde A\in L^\infty(X, T^*X)$. If 
$\|L^{1+\varepsilon}\tilde A\|_\infty\pp \infty$ (respectively
$\|L^{1+\varepsilon}V\|_\infty<\infty$) then the perturbation $\big(
d_{A_{\rm f}}^* (i\tilde A\wedge) + (i\tilde A\wedge)^* d_{A_{\rm
f}})$ (resp.\ $V$) is short-range in $\Bc(\Gr, \Gr^*)$. If $\tilde A$ is radial,
$\|L^{\varepsilon}\tilde A\|\pp \infty$ and
$\|L^{1+\varepsilon} x^{2-p}\partial_x \tilde A\|\pp \infty$, where
these norms are in $\Bc(L^2(X, g), L^2(X, \Lambda^1, g))$
(respectively $V$ radial, and $\|L^{1+\varepsilon} x^{2-p}\partial_x
V\|_\infty\pp \infty$) then the same perturbation is long-range in
$\Bc(\Gr, \Gr^*)$. 
\end{lemma} 
\proof
We deal with the magnetic perturbation. Start with the short-range. We
have $\langle L^{1+\varepsilon}(i\tilde A\wedge)^* 
d_{A_{\rm f}}) f\rangle= \langle d_A f, L^{1+\varepsilon}(i\tilde
A\wedge) f \rangle $ and on the other hand, we have
$\langle L^{1+\varepsilon} d_{A_{\rm f}}^* (i\tilde A\wedge) f,
f\rangle$ $=\langle [L^{1+\varepsilon}, d_{A_{\rm f}}^*]L^{-\varepsilon
} L^{\varepsilon} 
(i\tilde A\wedge) f, f\rangle + \langle L^{1+\varepsilon}(i\tilde
A\wedge) f, d_A f\rangle$. This is bounded by $\|f\| +\|d_A f\|^2$
uniformly in $f\in\cC^\infty_c(X)$.

We deal now with the long-range perturbation by checking
\eqref{e:lr'}. The condition with $L^\varepsilon$ is treated as above. 
In the variable of the free metric $g_p$ \eqref{gp'}, $\partial_r$ is
given by $\partial_L:=x^{2-p}\partial_x-(n-1)p/2\, 
x^{1-p}$. We extend
$\partial_L$ on $1$-forms by setting $\partial_L=x^{2-p}\partial_x-(n+1)p/2\, 
x^{1-p}$. Note it is symmetric on $1$-forms with compact support on the cusp. 
First, we have on smooth functions with compact support on the
cusp that:
\begin{equation}\label{e:comm1}
[d_{A_{\rm f}}, \partial_L]P_0= (2(1-p)x^{1-p}d_{A_{\rm f}} + c(1-p)x^{2(1-p)}x^{p-2}dx\wedge)P_0. 
\end{equation} 
Here, we used $d_{A_{\rm f}}P_0=(dx\wedge \partial_x\cdot)P_0$. Note that
$x^{p-2}dx\wedge$ is a bounded operator from function to $1$-forms. 
In the following, we drop $P_0$ and $\Xi$ to lighten the notation. For
$\varphi\in\cC^\infty_c(X)$, we have
\begin{eqnarray*}
\langle L^{1+\varepsilon} \varphi, [
d_{A_{\rm f}}^* (i\tilde A\wedge) + (i\tilde A\wedge)^* d_{A_{\rm
f}}] \varphi\rangle &=& \langle [\partial_L, d_{A_{\rm f}}]
L^{1+\varepsilon} \varphi, \tilde A\wedge \varphi\rangle + \\ 
&&\hspace*{-6cm}
+\langle \tilde A\wedge L^{1+\varepsilon} \varphi, [d_{A_{\rm f}},
\partial_L]\varphi\rangle + \langle d_{A_{\rm f}}
L^{1+\varepsilon}\varphi, [\tilde A\wedge, \partial_L]\rangle+
\langle[\partial_L, \tilde A\wedge] L^{1+\varepsilon}\varphi,
d_{A_{\rm f}}\varphi\rangle, 
\end{eqnarray*} 
Once $d_{A_{\rm f}}$ commuted with
$L^{1+\varepsilon}$, the two last terms are controlled by the
assumption on $\|[\partial_L, \tilde A\wedge]
L^{1+\varepsilon}\|$. The two first ones are $0$ for $p=1$ using
\eqref{e:comm1}. When $p\pp 1$, note that $x^{1-p}L^{1+\varepsilon }=c
L^{\varepsilon}$ and control the term using $L^\varepsilon \tilde A$
bounded. \qed 

We now describe the perturbation of the metric following the two classes. 
We keep the notation from Theorem \ref{t:mourre0}. We introduce the canonical 
unitary transformation due to the change of measure. Set
$\rho$ to be $\rho_{\rm sr}$, $\rho_{\rm lr}$ or $\rho_{\rm t}$. Let
$U$ be the operator of multiplication by $(1+\rho)^{-n/4}$ in $L^2(X,
g)$ and $V$ the operator of multiplication by $(1+\rho)^{(2-n)/4}$ in
$L^2(X,T^*X, g)$. The operator $U$ is a unitary operator from $ L^2(X,
g)$ onto $L^2(X, \tilde g)$ and $V$ a unitary operator from $L^2(X,
T^*X, g)$ onto $L^2(X, T^*X, \tilde g)$.

\begin{lemma}\label{l:uni}
Let $\widetilde\Delta_{A_{\rm f}}$ be the magnetic Laplacian of vector
potential $A_{\rm f}$ acting in $L^2(X, \tilde g)$. 
Let $W_0:= U^{-1}\widetilde\Delta_{A_{\rm f}} U$. Then,
\begin{enumerate}
\item On $\cC_c^{\infty}(X)$, $W_0$ acts by $U d_{A_{\rm f}}^* V^2
d_{A_{\rm f}} U$. In $L^2(X,g)$, it is essentially self-adjoint and
its domain is $\Dc(\Delta_{A_{\rm f}})$. 
\item One has $(W_0- \Delta_{A_{\rm f}} +i)^{-1}- (\Delta_{A_{\rm
f}}+i)^{-1}$ is compact. 
\item For $\rho_{\rm sr}$ and $\rho_{\rm t}$, $W_0$ is a short-range 
perturbation of the magnetic Laplacian $\Delta_{A_{\rm f}}$ in the space 
$\Bc(\Dc(\Delta_{A_{\rm f}}), \Dc(\Delta_{A_{\rm f}})^*)$.
\item For $\rho_{\rm lr}$, $W_0$ is a long-range perturbation of 
$\Delta_{A_{\rm f}}$ in $\Bc(\Dc(\Delta_{A_{\rm f}}),
\Dc(\Delta_{A_{\rm f}})^*)$. 
\end{enumerate} 
\end{lemma}
\proof
We write $\widetilde\Delta_{A_{\rm f}}$ with the help of the operator
$d_{A_{\rm f}}$. Since the manifold is complete,
$\widetilde\Delta_{A_{\rm f}}$ is essentially self-adjoint on
$\cC^{\infty}_c(X)$. In particular it corresponds to $\tilde
d_{A_{\rm f}}^*\tilde d_{A_{\rm f}}$, the Friedrichs extension.
This is equal to $U^2 d_{A_{\rm f}}^* V^2 d_{A_{\rm f}}$ in $L^2(X,
\tilde g)$. Now remark that $(1+\rho)^\alpha$ stabilizes
$\Dc(\Delta_{A_{\rm f}})$, $\Dc(d_{A_{\rm f}})$ and $\Dc(d_{A_{\rm
f}}^*)$, for all $\alpha\in\R$ to obtain the first point. 

We now compare the two operators in $L^2(X, g)$. We compute on
$\cC^\infty_c(X)$. 
\begin{equation}\label{e:D}\begin{split}
D:=W_0 - \Delta_{A_{\rm f}} &= U^{-1}
U^2 d_{A_{\rm f}}^* V^2 d_{A_{\rm f}} U - d_{A_{\rm f}}^* d_{A_{\rm f}}
\\
 &=
U^{-1}\widetilde\Delta_{A_{\rm f}}(U-1) +
U d_{A_{\rm f}}^* (V^2-1) d_{A_{\rm f}}+
(U-1)\Delta_{A_{\rm f}}. 
\end{split}\end{equation}
We focus on point $(3)$. The two first terms need a justification. We start with the first term.
\begin{equation*}
 L^{1+\varepsilon} U^{-1} \widetilde\Delta_{A_{\rm f}}(U-1)=
U^{-1} \big(\widetilde\Delta_{A_{\rm f}} + 
[L^{1+\varepsilon}, \widetilde\Delta_{A_{\rm f}}]L^{-1-\varepsilon}\big)
L^{1+\varepsilon}(U-1). 
\end{equation*} 
Using Lemma \ref{l:L} and the invariance of the domain under
$(1+\rho)^\alpha$, we obtain that $(\widetilde\Delta_{A_{\rm f}} +
[L^{1+\varepsilon}, \widetilde\Delta_{A_{\rm f}}]L^{-\varepsilon})^*$ is
bounded from $\Dc(\Delta_{A_{\rm f}})$ to $L^2(X, g)$. Again using properties
of $(1+\rho)^\alpha$, for all $\varphi\in\cC^\infty_c(X)$ we get
\begin{align*}
\lefteqn{|\langle \varphi, L^{1+ \varepsilon}
U^{-1}\widetilde\Delta_{A_{\rm f}}(U-1) \varphi \rangle |}\\&&&&=
|\langle (\widetilde\Delta_{A_{\rm f}} +
 [L^{1+\varepsilon}, \widetilde\Delta_{A_{\rm f}}]L^{-1-\varepsilon})^*
 U^{-1}\varphi, L^{1+\varepsilon }(U-1)\varphi \rangle |
\leq c \|(\Delta_{A_{\rm f}}+i) \varphi\|^2.
\end{align*} 
For the second term, we have $L^{1+ \varepsilon }\big(U
d_{A_{\rm f}}^* (V^2-1) d_{A_{\rm f}}\big)=U \big( d_{A_{\rm f}}^* +[L^{1+
 \varepsilon}, d_{A_{\rm f}}^*] L^{-1 - \varepsilon }\big)$ 
\linebreak
$L^{1+
 \varepsilon } (V^2-1) d_{A_{\rm f}}$. By Lemma \ref{l:L} and the
invariance of $\Dc(d_{A_{\rm f}})$ by $\rho^\alpha$, we obtain
\begin{align*}
\lefteqn{|\langle \varphi, L^{1+ \varepsilon }\big(U d_{A_{\rm f}}^* (V^2-1)
d_{A_{\rm f}}\big) 
\varphi\rangle|}\\&&&&= 
|\langle 
\big( d_{A_{\rm f}} + L^{-1 - \varepsilon } [d_{A_{\rm f}}, L^{1+
 \varepsilon}]\big)U \varphi, L^{1+
 \varepsilon } (V^2-1) d_{A_{\rm f}}\varphi\rangle|\leq c \|(\Delta_{A_{\rm f}}+i) \varphi\|^2, 
\end{align*} 
for all $\varphi\in\cC^\infty_c(X)$. To finish, use the fact that
$\cC^\infty_c(X)$ is a core for $\Delta_{A_{\rm f}}$. 

We now deal with $(4)$ by checking \eqref{e:lr'}. 
The real point to check is that $\|L^{1+\varepsilon }\tilde\Xi[D,
\partial_L]P_0\tilde\Xi\|_{\Bc(\Dc(\Delta), \Dc(\Delta)^*)}$ is
finite. We take $\partial_L$ like in the proof of Lemma \ref{l:magne}. First,
\begin{equation}\label{e:comm2}
[d^*d, \partial_L]P_0= c(1-p)x^{3(1-p)}P_0. 
\end{equation} 
We start with the easy part of $D$. We drop $\tilde\Xi$ and $P_0$ for
clarity. We have: 
\begin{align*}
\langle L^{1+\varepsilon }\varphi, [(U-1)d^*d, \partial_L] \varphi \rangle= 
 \langle L^{1+\varepsilon}\varphi, [U, \partial_L] d^*d \varphi
 \rangle+ \langle L^{1+\varepsilon}(U-1)\varphi, [ d^*d, \partial_L]
 \varphi \rangle. 
\end{align*} 
This is bounded by $\|(\Delta +i)\varphi\|^2$. Indeed, the first
term follows since $L^{1+\varepsilon}[U, \partial_L]$ is bounded in
$L^2(X,g)$. The second one is $0$ for $p=1$ and equals otherwise to $c
\langle \varphi, L^{\varepsilon -2}(U-1) \varphi\rangle$ by
\eqref{e:comm2}. Turn now to:
\begin{eqnarray*}
\lefteqn{\langle L^{1+\varepsilon}\varphi, [Ud^*(V^2-1)d,
\partial_L]\varphi\rangle}\\&
=&\langle dU L^{1+\varepsilon}\varphi, [V^2, \partial_L]d\varphi\rangle
+\langle [U, \partial_L] L^{1+\varepsilon}\varphi, d^*(V^2-1)d\varphi\rangle 
\\&&
+\langle (V^2-1)dU L^{1+\varepsilon}\varphi, [d, \partial_L]
\varphi\rangle + \langle [d, \partial_L] U L^{1+\varepsilon}\varphi,
(V^2-1)d\varphi \rangle. 
\end{eqnarray*} 
The first is controlled by commuting $L^{1+\varepsilon }$ with $d$ like
above and by using that $L^{1+\varepsilon }[V^2, \partial_L]$ is
bounded in $L^2(X, \Lambda^1, g)$. For the second one, 
$[U, \partial_L] L^{1+\varepsilon}$
is boun\-ded in $L^2(X, \Lambda^1, g)$. Turn the two last ones
and use \eqref{e:comm2}, for $p=1$ this is $0$. Focus on the very last
one for example. Now commute $L^{1+\varepsilon}$ with $d$ like
above. The most singular term being $\langle d U \varphi,
x^{1-p}L^{1+\varepsilon}(V^2-1) d\varphi\rangle$. Now remember that
$x^{1-p}L^{1+\varepsilon }=c L^{\varepsilon}$ and use 
the fact that 
$L^\varepsilon (V^2-1)$ is bounded to control it. To conclude, repeat the
same arguments for $[U^{-1}\tilde\Delta(U-1), \partial_L]$. 

We turn to point (2), $W_0$ and $\Delta_{A_{\rm f}}$ have the same
domain. We take the proof of Lemma \ref{l:pertuA} replacing $\cG$ with this
domain. We then obtain a rigorous version of \eqref{33}. Therefore, it remains
to check that $W_0-\Delta_{A_{\rm f}}\in\Kc(\Dc(\Delta_{A_{\rm f}}),
\Dc(\Delta_{A_{\rm f}})^*)$. This comes directly using Rellich-Kondrakov Lemma
and \eqref{e:D}. \qed 

Finally, we gather various technicalities concerning the operator $L$. 

\begin{lemma}\label{l:L}
We have 
that $dL$ is with support in $(0, \varepsilon)\times M_{\rm nt}$ and $dL=
f(x)dx$ where $f:(0, \varepsilon )\rightarrow \R$ such that $f$ is $0$ in a
neighborhood of $\varepsilon$ and such that $f(x)= -x^{p-2}$ for $x$ small enough. Moreover:
\begin{enumerate}
\item The operator $L^{-\varepsilon }d(L^{1+\varepsilon })\wedge$ belongs to
 $\Bc\big(L^2(X, g), L^2(X,T^*X, g)\big)$ and the commutator
 $L^{-\varepsilon }[\Delta_{A_{\rm f}}, L^{1+
 \varepsilon}]$ with initial domain $\cC^\infty_c(X)$ extends to 
a bounded operator in 
 $\Bc\big( \Dc(\Delta_{A_{\rm f}}), L^2(X, g)\big)$.
\item $e^{itL}\Dc(\Delta_{A_{\rm f}})\subset \Dc(\Delta_{A_{\rm f}})$
and $\|e^{itL}\|_{\Bc(\Dc(\Delta_{A_{\rm f}}))}\leq c 
(1+t^2)$. 
\item $L^{-1-\varepsilon }\Dc(\Delta_{A_{\rm f}})\subset
\Dc(\Delta_{A_{\rm f}})$. 
\end{enumerate} 
\end{lemma} 
\proof
With the diagonalization of Section \ref{s:diago}, the operator 
$\Delta_{A_{\rm f}}$ is given by \eqref{e:H_0}.
The operator $L$ corresponds to the operator $L_0$ of multiplication by $r\otimes 1_{M_{\rm nt}}$ on $(c, \infty)$ in this variable and by $1$ on the rest of the manifold. 
Hence, points (1) and (3) are easily obtained. 
Moreover $e^{itL_0}/(1+t^2)$ and its first and second derivative belong to 
$L^2(X_0)$, uniformly in $t$, from which (2) follows.\qed

\section{The non-trapping case for perturbed metrics} \label{s:ntp}
\subsection{The essential spectrum}

In this section, we compute the essential spectrum of a magnetic
Laplacian given by a non-trapping vector potential. Unlike the
trapping case, it is non-empty in the complete case. To show this, we 
apply perturbation techniques to the results of the previous section.

We restrict ourselves to conformal perturbations of 
\emph{exact} cusp metrics. In a fixed product decomposition of $X$ 
near $M$ we rewrite \eqref{cume} as 
\begin{equation}\label{cume2} 
g_0=a \left(\frac{dx}{x^2}+\alpha(x)\right)^2+h(x) 
\end{equation} 
where $a\in\cun(\oX)$, $\alpha\in\cun([0,\varepsilon)\times M,\Lambda^1(M))$
and $h\in \cun([0,\varepsilon)\times M, S^2TM)$. 
By \cite[Lemma 6]{wlom}, the function $a_0:=a(0)$, 
the metric $h_0:=h(0)$ and the class (modulo exact forms) 
of the $1$-form $\alpha_0:=\alpha(0)$, defined on $M$, 
are independent of the chosen product decomposition and of 
the boundary-defining function $x$ inside the fixed cusp structure. 
 
\begin{definition}\label{defex} 
The metric $g_0$ is called \emph{exact} if $a_0=1$ and $\alpha_0$ is an 
exact $1$-form. 
\end{definition} 
 
If $\alpha_0=df$ is exact, then by replacing $x$ with the boundary-defining 
function $x'=x/(1+xf)$ inside the same cusp structure, we can as 
well assume that $\alpha_0=0$ (see \cite{wlom}). It follows that 
$g_0$ is quasi-isometric to a cylindrical metric near infinity. 
 
\begin{proposition}\label{p:thema} 
Let $(X,\tilde g_p)$ be a Riemannian manifold 
with a conformal \emph{exact} cusp metric $\tilde g_p:=(1+\rho)g_p$, where
$g_0$ is exact, $g_p=x^{2p}g_0$ and 
\begin{equation*}
\rho\in L^\infty(X; \R), \quad \inf_{x\in X}\rho(x)\pg -1, \quad 
\rho(x)\rightarrow 0,
\mbox{ as } x\rightarrow 0.
\end{equation*} 
Let $A$ be a non-trapping vector potential given by \eqref{aft}. Then 
\begin{itemize} 
\item For $0<p\leq 1$, the Friedrichs extension of
$\Delta_A$ has essential spectrum $ \siges(\Delta_A)= [\kappa(p), \infty)$, 
where $\kappa(p)=0$ for $p<1$ and $\kappa(1)=(n-1)/2$. Moreover, if $\rho$
is smooth, the $\Delta_A$ is essential self-adjoint on $\cC^\infty_c(X)$.
\item If $p>1$ and $\tilde g_p:=g_p$ is the unperturbed metric given
in \eqref{mc}, then every self-adjoint extension of $\Delta_A$
has empty essential spectrum. 
\end{itemize} 
\end{proposition} 
Note that for $p>1$, the unperturbed metric $g_p$ given in \eqref{mc} 
is essentially of \emph{metric horn} type \cite{lepe}.

\proof
Using the Weyl theorem, Lemma \ref{l:pertuA} and by changing the
gauge, we can suppose without loss of generality that $A$ is of
the form \eqref{acon} with $C=0$.

We start with the complete case. The essential self-adjointness follows
from \cite{shubin}. In the exact case, $g_p$ is quasi-isometric to 
the metric \eqref{gp'}. Using \cite[Theorem 9.4]{GG} (see Theorem 9.5
for the case of the Laplacian), to compute
the essential spectrum we may replace $h(x)$ in
\eqref{cume2} by the metric $h_0:=h(0)$ on $M$, extended to a
symmetric $2$-tensor constant in $x$ near $M$, and we may set $\rho=0$. 
By Proposition \ref{p:cut}, computing
$\siges(\Delta_A)$ is the same as computing $\siges(\Delta^l_A)$ on
$X'$ of Proposition \ref{p:free'}. By the results of Section
\ref{s:diago}, this is given on $L^2(\R^+)$ by $\siges(-\Delta+V_p)
=[V_p(\infty), \infty)$. 

Let now $p> 1$. The metric is no longer complete and $(X,g_p)$ is not proper;
one can not apply \cite[Theorem 9.4]{GG}. By Lemma \ref{l:indices} and
by the Krein formula, all self-adjoint extensions have the same essential
spectrum. So it is enough to consider the Friedrichs extension of
$\Delta_A$. We now use Propositions \ref{p:free'} and
\ref{p:cut}. The operator $D^*D$ is non-negative, so the spectrum of
$\Delta_A^0$ is contained in $[\varepsilon^{2-2p}c_0,\infty)$. By
Proposition \ref{p:cut}, the essential spectrum does not depend on
the choice of $\varepsilon$. Now we remark that $p>1$ implies
$c_0\neq 0$. Indeed, the equality would imply that $1/p\in\zz$, which
is impossible. Thus by letting $\varepsilon \to 0$ we conclude that the
essential spectrum is empty. \qed 

We have used above a general lemma about compact perturbations of magnetic 
Laplacians:

\begin{lemma}\label{l:pertuA}
Let $A$ and $A'$ in $L^\infty(X, T^*X)$ be two magnetic fields on a
smooth Riemannian manifold $(X,g)$ (possibly incomplete) with a measurable 
metric $g$. Suppose that $A-A'$ belongs to $L^\infty_0(X,
T^*X)$. Let $\Delta_A={d_A}^*d_A$ and
$\Delta_{A'}={d_{A'}}^*d_{A'}$. Then $(\Delta_A+i)^{-1}- (\Delta_{A'}
+i)^{-1}$ is compact. 
\end{lemma} 
Here, $L^\infty_0(X, T^*X)$ denotes the space of those forms of 
$L^\infty(X,T^*X)$ which are norm limit of compactly supported forms. Note
that this lemma holds without any modification for 
a $C^1$ manifold equipped with a (RM) structure, see \cite[section
9.3]{GG}. Unlike the result on the stability of the essential
spectrum of the (magnetic) Laplacian from \cite[Theorem 9.5]{GG}, we do
not ask for the completeness of the manifold. A magnetic perturbation
is much less singular than a perturbation of the metric. 
\proof 
Note that the form domain of $\Delta_A$ and of $\Delta_A$ is given by
$\Gr:=\Dc(d)$, because $A$ and $A'$ are in $L^\infty(X,
T^*X)$. We write with a tilde the extension of the magnetic Laplacians to 
$\Bc(\Gr, \Gr^*)$. We aim to give a rigorous meaning to 
\begin{equation}\label{33}
(\Delta_A+i)^{-1}-(\Delta_{A'}+i)^{-1}=(\Delta_A+i)^{-1}(\Delta_{A'}-\Delta_A) (\Delta_{A'}+i)^{-1}.
\end{equation}
We have $(\Delta_A+i)^{-1*}\Hr\subset\Gr$. This allows one to
deduce that ($\Delta_A+i)^{-1}$ extends to a unique continuous operator
$\Gr^*\rightarrow\Hr$. We denote it for the moment by $R$. From
$R(\Delta_A+i)u=u$ for $u\in \Dc(\Delta_A)$ we get, by density of 
$\Dc(\Delta_A)$ in $\Gr$ and continuity, $R(\widetilde{\Delta_A}+i)u=u$ for
$u\in\Gr$, in particular
\begin{equation*}
(\Delta_{A'}+i)^{-1}= R(\widetilde{\Delta_A}+i)(\Delta_{A'}+i)^{-1}.
\end{equation*}
Clearly,
\begin{equation*}
(\Delta_A+i)^{-1}= (\Delta_A+i)^{-1}(\Delta_{A'}+i)(\Delta_{A'}+i)^{-1}=
R(\widetilde{\Delta_{A'}}+i)(\Delta_{A'}+i)^{-1}.
\end{equation*}
We subtract the last two relations to get
\begin{equation*}
(\Delta_A+i)^{-1}-(\Delta_{A'}+i)^{-1}=R(\widetilde{\Delta_{A'}}
-\widetilde{\Delta_A})(\Delta_{A'}+i)^{-1} 
\end{equation*}
Since $R$ is uniquely determined as the extension of $(\Delta_A+i)^{-1}$ to a
continuous map $\Gr^*\rightarrow\Hr$, one may keep the notation
$(\Delta_A+i)^{-1}$ for it. With this convention, the rigorous version of
(\ref{33}) that we shall use is:
\begin{equation}\label{34}
(\Delta_A+i)^{-1}-(\Delta_{A'}+i)^{-1}=(\Delta_A+i)^{-1}(\widetilde
\Delta_{A'}-\widetilde \Delta_A)(\Delta_{A'}+i)^{-1}. 
\end{equation}
Since $A''=i( A-A')\in L^\infty_0(X, T^*X)$, the Rellich-Kondrakov
lemma gives that $A''\wedge$ belongs to $\Kc(\Gr, L^2(X, T^*X))$. Therefore,
$\widetilde\Delta_{A'}-\widetilde\Delta_{A} = (A''\wedge)^* d_{A} - 
d_{A}^*A''\wedge\in \Kc(\Gr,\Gr^*)$. This gives the announced
compactness.\qed 

\subsection{The spectral and scattering theory}
In this section, we refine the study of the essential spectrum given
in Proposition \ref{p:thema} for non-trapping vector potential.
As the essential spectrum arises only in the complete case, we will suppose
that $p\leq 1$. 

We give below our main result in the study of the nature of the essential
spectrum and in scattering theory under short-range perturbation. It 
is a consequence of the Mourre theory \cite{mou} with an
improvement for the regularity of the boundary value of the resolvent,
see \cite{ggm} and references therein. 

We treat some conformal perturbation of the metric $\eqref{gp'}$. To our
knowledge, this is the weakest hypothesis of perturbation of a metric 
obtained so far using Mourre theory. Compared to previous approaches, we use a
conjugate operator which is local in energy and therefore can be compared
directly to a multiplication operator. We believe that this procedure can
be implemented to all known Mourre estimates on manifolds to improve the
results obtained by perturbation of the metric. 

We fix $A_{\rm f}$ a non-trapping vector potential of the form \eqref{acon}. 
By a change of gauge, one can suppose that $A_{\rm f}=\theta_0$ which 
is constant in a neighborhood of $M$. Let $M_{\rm t}$ (resp.\ $M_{\rm nt}$) 
be the union of the connected components of $M$ on which $A_{\rm f}$ 
is trapping (resp.\ non-trapping). Let $L$ be the operator 
of multiplication by a smooth function $L\geq 1$ which measures the length 
of a geodesic going to infinity in the directions where $\theta_0$ is 
non-trapping: 
\begin{eqnarray}\label{e:L} 
 L(x)=\begin{cases} 
-\ln(x)& \text{for $p=1$},\\ 
-\frac{x^{p-1}}{p-1}& \text{for $p< 1$} 
\end{cases}
\end{eqnarray} 
on $(0, \varepsilon/4)\times M_{\rm nt}$ for small $x$, 
and $L=1$ on the trapping part
$\big((0, \varepsilon /2)\times M_{\rm nt} \big)^c$. Given $s\geq 0$, 
let $\Lr_{s}$ be the domain of $L^s$ equipped with the graph norm. We set
$\Lr_{-s}:=\Lr_{s}^*$. Using the Riesz theorem, we obtain the scale of
spaces $\Lr_{s}\subset L^2(X, \tilde g) \subset \Lr_{s}^*$, 
with dense embeddings and where $\tilde g$ is defined in the theorem below.
Given a subset $I$ of 
$\R$, let $I_{\pm}$ be the set of complex number $x\pm iy$, where
$x\in I$ and $y>0$. The thresholds $\{\kappa(p)\}$ are given in Proposition
\ref{p:thema}.

For shorthand, perturbations of \emph{short-range} type (resp.\
\emph{trapping} type) are denoted with the subscript ${\rm sr}$
(resp.\ ${\rm t}$); they are supported in $(0, \varepsilon)\times
M_{\rm nt}$ (resp.\ in $((0, \varepsilon/2)\times M_{\rm nt})^c$). 
We stress that the class of ``trapping type''
perturbations is also of short-range nature, in the sense
described in Section \ref{s:srlr}, even if no decay is required. This
is a rather amusing phenomenon, linked to the fact that no essential
spectrum arises from the trapping cusps. The subscript $\rm lr$
denotes \emph{long-range} type perturbations, also
with support in $(0, \varepsilon)\times M_{\rm nt}$. We ask such 
perturbations to
be \emph{radial}, i.e.\ independent of the variables in $M$. In other words, a
perturbation $W_{\rm lr}$ satisfies $W_{\rm lr}(x,m)= W_{\rm
lr}(x,m')$ for all $m,m'\in M$. 

\begin{theorem}\label{t:mourre0}
Fix $\varepsilon >0$. Let $X$ be endowed with the metric $\tilde g =
(1+\rho_{\rm sr}+\rho_{\rm lr}+\rho_{\rm t}) g$, where $g=g_p$ is given in
\eqref{gp'} for some $0<p\leq 1$; $\rho_{\rm sr}$, $\rho_{\rm lr}$ and 
$\rho_{\rm t}$ belong to $\cC^\infty(X)$, such that 
\begin{align*}
\inf_{x\in X}(\rho_{\rm sr}(x)+ \rho_{\rm lr}(x) + \rho_{\rm
t}(x))\pg -1,&& \rho_{\rm t}(x) = o(1) \text{ as }
x\rightarrow 0\end{align*}
 and such that
\begin{align*}
L^{1+\varepsilon}\rho_{\rm sr}, 
d_{A_{\rm f}}\rho_{\rm sr}, 
\Delta_{A_{\rm f}}\rho_{\rm sr}, 
L^{\varepsilon}\rho_{\rm lr},
L^{1+\varepsilon}d_{A_{\rm f}}\rho_{\rm lr},
\Delta_{A_{\rm f}}\rho_{\rm lr},
d_{A_{\rm f}}\rho_{\rm t},
\Delta_{A_{\rm f}}\rho_{\rm t} \text{ belong to $L^\infty$.}
\end{align*} 
In $\Hr= L^2(X, \tilde g)$, let $\widetilde\Delta_A$ be the magnetic
Laplacian with a non trapping potential 
$A= A_{\rm f} + A_{\rm lr } + A_{\rm sr}+ A_{\rm t}$, where 
$A_{\rm f}$ is as in \eqref{acon}, 
$A_{\rm sr}$, $A_{\rm lr}$ and $A_{\rm t}$ 
are in $\cC^\infty(X, T^*X)$ such that: 
\begin{align*}
 \|L^{1+\varepsilon}A_{\rm sr}\|_\infty,\,
 \|L^{\varepsilon}A_{\rm lr}\|_\infty,\,
 \|L^{1+\varepsilon}\cL_{x^{(2-p)}\partial_{x}} A_{\rm lr}\|_\infty\pp 
\infty \mbox{ and } A_{\rm t}= o(1), 
\end{align*} 
where $\cL$ denotes the Lie derivative.
Let $V= V_{\rm loc}+ V_{\rm sr}+ V_{\rm lr} + V_{\rm t}$ and $V_{\rm
lr}$ be some potentials, where $V_{\rm loc}$ is measurable with
compact support and $\widetilde\Delta_{A_{\rm f}}$-compact and
$V_{\rm sr}$, $V_{\rm lr}$ and $V_{\rm t}$ are in $L^\infty(X)$ such
that: 
\begin{equation*}
\|L^{1+\varepsilon}V_{\rm sr}\|_\infty,\,
 \|L^{1+\varepsilon}d_{A_{\rm f}}V_{\rm lr}\|_\infty<\infty 
\mbox{ and } V_{\rm t}=o(1). 
\end{equation*} 
Consider the magnetic Schr\"odinger operators $H_0=
\widetilde\Delta_{A_{\rm f}+A_{\rm lr}}+V_{\rm lr}$ and $H=
\widetilde\Delta_{A} + V$. Then 
\begin{enumerate}
\item $H$ has no singular continuous spectrum.
\item The eigenvalues of $\R\setminus \{\kappa(p)\}$ have finite
multiplicity and no accumulation points outside $\{\kappa(p)\}$.
\item Let $\cJ$ a compact interval such that $\cJ\cap \big(\{\kappa(p)\}\cup
 \sigma_{\rm pp}(H)\big)=\emptyset$. Then, for all $s\in (1/2, 3/2)$,
there exists $c$ such that 
\begin{equation*}
\|(H-z_1)^{-1} - (H-z_2)^{-1} \|_{\Bc(\Lr_s, \Lr_{-s})} \leq c
\|z_1-z_2\|^{s-1/2}, 
\end{equation*}
for all $z_1, z_2 \in \cJ_{\pm}$.
\item Let $\cJ=\R\setminus \{\kappa(p)\}$ and let $E_0$ and $E$ be the
continuous spectral component of $H_0$ and $H$, respectively. Then,
the wave operators defined as the strong limit
\begin{equation*}
\Omega_{\pm}=\slim_{t\rightarrow \pm \infty} e^{itH}e^{-itH_0}E_0(\cJ) 
\end{equation*} 
exist and are complete, i.e.\ $\Omega_{\pm} \Hr= E(\cJ)\Hr$.
\end{enumerate} 
\end{theorem}

\begin{remark}
Any smooth $1$-form $A$ on $\oX$
is a short-range perturbation of a free vector 
potential $A_{\rm f}$ as in \eqref{acon}.
\end{remark}

\begin{remark}
If one is interested only in the free metric $g_p$, 
the conclusions of the theorem hold for $V_{\rm loc}$ in
the wider class of $\widetilde\Delta_{A_{\rm f}}$-form compact perturbations,
using \cite[Theorem 7.5.4]{ABG}. Similar results should hold
for a smooth metric on $\oX$ and Dirichlet boundary conditions, 
as in \cite{bouclet}.
\end{remark}

\proof 
We start with $\Delta_{A_{\rm f}}$ in $L^2(X, g)$. In section \ref{s:conjop},
we transform it unitarily into $\Delta_0$ given by \eqref{e:H_0}. For
$R$ finite, we construct a conjugate operator $S_R$ to $\Delta_0$ given by
\eqref{e:SR}. Theorem \ref{t:mourre_0} gives a Mourre estimate for
$\Delta_0$ and the regularity of $\Delta_0$ compared to the self-adjoint
operator $S_R$. We go back by unitary transform into $L^2(X,
g)$. Since the dependence on $R$ is no longer important, we denote
simply by $S$ the image of the conjugate operator $S_R$. Therefore, we
have $\Delta_{A_{\rm f}}\in \cC^{2}(S, \Dc(\Delta_{A_{\rm f}}), L^2(X,
g))$ and given $\cJ$ an open interval included in $\sigma_{\rm
ess}(H_0)$, there is $c>0$ and a compact operator $K$ such that 
the inequality 
\begin{equation}\label{e:T}
E_\cJ(T) [T,iS] E_\cJ(T) \geq c E_\cJ(T)+ K 
\end{equation}
holds in the sense of forms in $L^2(X,g)$, for $T=\Delta_{A_{\rm f}}$.

Let $W_0$ be the unitary conjugate of $\widetilde\Delta_{A_{\rm f}}$ 
acting in $L^2(X,g)$. By Lemma \ref{l:uni}, 
$W_0\in\cC^{1,1}(S,\Dc(\Delta_{A_{\rm f}}), \Dc(\Delta_{A_{\rm f}})^*)$
is a sum of short and long-range perturbation as described in 
Section \ref{s:srlr}. In
particular, we get $W_0\in\cC^{1}_{\rm u}(S,\Dc(\Delta_{A_{\rm f}}),
\Dc(\Delta_{A_{\rm f}})^*)$. By the point (2) of Lemma \ref{l:uni} and
\cite[Theorem 7.2.9]{ABG} the inequality \eqref{e:T} holds for 
$T$ (up to changing $c$ and $K$). 

We now go into $L^2(X, \tilde g)$ using $U$ defined before Lemma
\ref{l:uni}. We write the con\-ju\-ga\-te operator obtained in this way by
$\tilde S$. Therefore, $\tilde\Delta_{A_{\rm f}}$ belongs to $\cC^{1,1}(\tilde
S, \Dc(\Delta_{\tilde A_{\rm f}}), \Dc(\Delta_{\tilde A_{\rm f}})^*)$ and
given $\cJ$ an open interval included in $\sigma_{\rm ess}(H_0)$,
there is $c>0$ and a compact operator $K$ such that 
\begin{equation}\label{e:tildeT}
E_\cJ(\tilde T) [\tilde T,i\tilde S] E_\cJ(\tilde T) \geq c E_\cJ(\tilde T)+ K 
\end{equation}
holds in the sense of forms in $L^2(X, \tilde g)$ for $\tilde T=\tilde
\Delta_{A_{\rm f}}$. We now add the perturbation given by 
$A_{\rm sr}, A_{\rm lr}, A_{\rm t}, V_{\rm sr}, V_{\rm lr},$ and $V_{\rm t}$. 
Note that $H$ has the same domain as $H_0$ and that $(H+i)^{-1}- (H_0+i)^{-1}$ 
is compact by Rellich-Kondrakov lemma and Lemma \ref{l:pertuA}. By Lemma
\ref{l:magne}, we obtain $H\in\cC^{1,1}(S,\Dc(H),\Dc(H)^*)$. As above, 
the inequality \ref{e:tildeT} is true for $\tilde T=H$. 

We now deduce the different claims of the theorem. The first comes from
\cite[Theorem 7.5.2]{ABG}. The second ones is a consequence of the Virial
theorem. For the third point first note that $\Lr_s\subset\Dc(|A|^s)$ for
$s\in [0,2]$ by Lemma \ref{l:sa} and use \cite{ggm} for instance (see
references therein). Finally, the last point follows from \cite[Theorem
 7.6.11]{ABG}. 
\qed

\section{The non-stability of the essential spectrum and of the
situation of limiting absorption principle}\label{s:nonstab} 

In $\rz^n$, with the flat metric, it is well-known that
only the behavior of the magnetic field at infinity plays
a r\^ole in the computation of the essential spectrum.
Moreover, \cite[Theorem 4.1]{KS} states that the non-emptiness
of the essential spectrum is preserved by the addition of a bounded magnetic
field, even if it can become purely punctual. Concerning a compactly
support magnetic field, the essential spectrum remains the same, see
\cite{MS}. However, it is well-known that one obtains a long-range
effect from it, in other words it acts on particles which have support
away from it. In the case of $\rz^n$ with a hole of some kind, this
phenomena are of special physical interests and are related to the
Aharonov-Bohm effect, see section \ref{s:AB} and references therein.

In contrast with the Euclidean setting, Theorem \ref{t:thmag}
indicates that in general the essential spectrum may vanish under
\emph{compactly supported} perturbations of the magnetic field. 
In the next sections, we discuss this effect both with and without 
the hypothesis of gauge invariance, and we investigate the coupling 
constant effect. 

\subsection{The case $H^1(X)= 0$.}\label{s:inv} 

In this section, we assume gauge invariance. We first characterize
trapping condition in terms of the magnetic field, see
Definition \ref{def7} for the case of a magnetic potential. 
We recall that if $H^1(X)= 0$, given a magnetic field $B$ the spectral
properties of the magnetic Laplacian $\Delta_A$ will not depend on the
choice of vector potential $A$ such that $dA=B$. Indeed, given $A,
A'$ such that $dA=dA'$, the operators $\Delta_A$ and $\Delta_{A'}$
are unitarily equivalent by a gauge transformation. Therefore, we denote the
magnetic Laplacian by $\Delta_B$ 
and express the condition of being (non-)trapping in function of $B$.

Let $p>0$ and let $X$ be the interior a compact manifold $\overline{X}$
 endowed with the metric $g_p$ given by
\eqref{pme}. For simplicity, assume that $B$ is a smooth $2$-form on
$\oX$ such that its restriction to $X$ is exact. Then there exists $A\in
\cun(\oX,T^*\oX)$ such that $B=dA$ (since the cohomology of
the de Rham complex on $\oX$ equals the singular cohomology of $\oX$,
hence that of $X$). Let 
\begin{align*}M=\sqcup_{\alpha\in\cA}M_\alpha
\end{align*} 
be the decomposition of the boundary $M$ into its connected components. 
Set 
\[\cA_0:=\{\alpha\in\cA; H^1(M_\alpha;\rz)=0\}.\]
For some $\cB\subset \cA$ set $M_\cB=\sqcup_{\beta\in\cB}M_\beta$ and 
consider the long exact cohomology sequence of the pair
$(\oX,M_\cB)$ with real coefficients:
\[H^1(\oX;\rz)\longrightarrow H^1(M_\cB;\rz)\stackrel{\partial}{\longrightarrow} 
H^2(\oX,M_\cB;\rz)\stackrel{i}{\longrightarrow} H^2(\oX;\rz)\]
Since we assume that $H^1(X;\rz)=0$ it follows that the connecting map
$\partial$ is injective. If $B$ vanishes under pull-back to $M_\cB$ then 
(since it is exact on $X$) it defines a class in $H^2(\oX,M_\cB;\rz)$ 
which vanishes under the map $i$, so it belongs to the image of the
injection $\partial$. We denote by $[B]_\beta$ the component of $[B]$ 
inside $\partial H^1(M_\beta)\subset H^2(\oX,M_\cB;\rz)$.

\begin{definition}\label{magtrap}
Assume $H^1(X)=0$. Let $B$ be a smooth exact $2$-form on $\oX$. Denote
by $\cB$ the set of those $\alpha\in\cA$ such that $B$ vanishes 
identically on $M_\alpha$. The field $B$ is called \emph{trapping} if
for each $\beta\in\cB$, 
the component $[B]_\beta \in\partial H^1(M_\beta)\subset H^2(\oX,M_\cB;\rz)$ 
is not integral, i.e.,
it does not live in the image of the map of multiplication by $2\pi$
\[H^2(\oX,M_\cB;\zz)\stackrel{2\pi\cdot}{\longrightarrow}H^2(\oX,M_\cB;\rz),\]
and \emph{non-trapping} otherwise.
\end{definition}

This definition is consistent with Definition \ref{d:intro} 
when $M$ is connected. Note that if $B$ is trapping then
$\cB$ must contain the index set $\cA_0$ defined above.

In order to apply Theorem \ref{t:mourre0} and Theorem \ref{t:thmag} 
we use the following lemma:

\begin{lemma}\label{l:compsupp}
\begin{enumerate}
\item \label{compsupp1} Let $A$ be a smooth vector potential on $\oX$ 
such that $dA=0$ in a neighbourhood of $M=\partial X$. Then
there exists a smooth vector potential $A'$, \emph{constant in $x$} in a 
neighborhood of $M$, such that $A=A'$ on $M$ and 
$d(A-A')=0$.
\item \label{compsupp2} Assume that $H^1(\oX,\rz)$ vanishes. 
Let $B$ be a trapping magnetic field on $\oX$. Then every vector potential
for $B$ will be trapping.
\item \label{compsupp3} Assume moreover that $H_1(X;\zz)=0$. 
Let $B$ be non-trapping such that $\imath^*_M B=0$. 
Then every vector potential
for $B$ will be non-trapping.
\end{enumerate}
\end{lemma}

Recall that $\pi_1(X)=0\Longrightarrow H_1(X;\zz)=0 
\Longrightarrow H^1_\dR(X)=0$.

\begin{proof}
\ref{compsupp1}) Let us show that one can choose $A$ to be
\emph{constant in $x$} 
near the boundary, in the sense that near $M$ it is the pull-back of a
form from $M$ under the projection $\pi:[0,\varepsilon)\times M\to M$
for $\varepsilon$ small enough. Indeed, 
$A-\pi^* \imath_M^*A$ is closed on the cylinder $[0,\varepsilon)\times M$ 
and vanishes when pulled-back to $M$. Now $M$ is a deformation-retract 
of the above cylinder, so the 
map of restriction to $M$ induces an isomorphism in cohomology and thus the 
cohomology class 
$[A-\pi^* \imath_M^*A]\in H^1([0,\varepsilon)\times M)$ must be zero.
Let $f\in\cun(\oX,\rz)$ be a primitive of this form for $x\leq \varepsilon/2$,
then $A-df$ is the desired constant representative.

\ref{compsupp2}) Consider the commutative diagram 
\begin{equation}\label{cohoz}
\begin{CD}
H^{1}(M_\cB;\rz)@>{\partial}>>H^{2}(\oX,M_\cB;\rz)@>>>H^2(\oX,\rz)\\
 @AA{2\pi\cdot}A@AA{2\pi\cdot}A@AA{2\pi\cdot}A\\
H^{1}(M_\cB;\zz)@>{\partial}>>H^{2}(\oX,M_\cB;\zz)@>>>H^2(\oX,\zz)
\end{CD}
\end{equation}
where the horizontal maps come form the long exact sequence of the
pair $(\oX,M)$ and the vertical maps are multiplication by $2\pi$.
Let $B$ be trapping and choose a vector potential $A$ for $B$, smooth on $\oX$.
We claim that $A$ is trapping.
Indeed, $A$ is not closed on the components $M_\alpha$, 
$\alpha\in\cA\setminus\cB$, while it is closed on $M_\cB$. We note that
$\partial[A_{|M_\cB}]=[dA]=[B]$, so $\partial [A_{|M_\beta}]=[B]_\beta$. 
Assume that for some $\beta\in\cB$, the class $[A_{|M_\beta}]$ were integral. 
Then using the first square from diagram \eqref{cohoz}, it would follow that 
$[B]_\beta$ was also integral, contradiction.

\ref{compsupp3}) If
$B$ is non-trapping, there exists $\beta\in\cB$ and $b\in H^2(\oX,M_\cB,\zz)$
with $[B]_\beta=2\pi b$.
The image of $[B]_\beta\in\partial H^1(M_\beta,\rz)$ in $H^2(\oX,\rz)$
is zero, thus $b$ maps to a torsion element in $H^2(\oX,\zz)$.
From $H_1(X,\zz)=0$ we see using the universal coefficients 
theorem 
\[0\to \Ext(H_1(\oX,\zz),\zz)\to H^2(\oX;\zz)\to \Hom(H_2(\oX,\zz),\zz)\to 0\]
that $H^2(\oX;\zz)$ is torsion-free.
Thus $b$ comes from some $a\in H^1(M_\cB,\zz)$. By commutativity we have 
$\partial(2\pi a)=[B]_\beta$. Since $\partial$ is an injection, it follows that
every vector potential smooth on $\oX$ for $B$ will define an integral 
cohomology class on $M_\beta$, thus will be non-trapping.
\end{proof}

A spectacular example is
a compactly supported magnetic field which induces very
strong long-range effects. If $H^1(\oX;\zz)=0$ and 
$B$ be is an exact $2$-form with compact support
in $X$, then $B$ is maximal non-trapping (see Definition \ref{def7}) 
if and only if its class in
$H^2(\oX,M;\rz)$ is integral. We will stress more on this aspect in the
following section.

We can also construct
compactly supported magnetic fields for which the consequences of
Theorem \ref{t:thmag} hold true. We
summarize this fact in the next proposition and give an explicit
construction in the proof. 

\begin{proposition}\label{p:contreex}
Let $X$ be the interior of a compact manifold $\oX$ with
boundary $M=\partial \oX$, endowed with a conformally cusp metric $g_p$. 
Assume that
$H^1(X)=0$ and $H^1(M_j)\neq 0$ for every connected component $M_j$ of
the boundary. Then there exists a non-zero smooth magnetic field $B$ with
\emph{compact support} such that the essential spectrum of $\Delta_B$ is
empty and such that for $p\geq 1/n$ the growth law of the eigenvalues does not
depend on $B$ and is given by \eqref{e:thmag}. Such fields $B$ are generic 
inside compactly supported magnetic fields.
\end{proposition}

\proof We construct $A$ like in (\ref{aft}) satisfying the hypotheses
of Theorem \ref{t:thmag}. We take $\varphi_0$ to be constant. Let
$\psi\in\cun([0, \varepsilon))$ be a cut-off function
such that $\psi(x)=0$ for $x\in[3\varepsilon/4,\varepsilon)$ and
$\psi(x)=1$ for $x\in[0, \varepsilon/2)$. Since $H^1(M_j)\neq 0$,
there exists a closed $1$-form $\beta_j$ on $M$ which is not exact.
Up to multiplying $\beta_j$ by a real constant, we can assume that the
cohomology class $[\beta_j]\in H^1_{\mathrm{dR}}(M_j)$ does not belong
to the image of $2\pi H^1(M_j;\zz)\to H^1(M_j;\rz)\simeq
H^1_{\mathrm{dR}}(M_j)$. Let $\beta$ denote the form on $M$ which equals
$\beta_j$ on $M_j$. Choose $A$ to be $\psi(x)\beta$ for
$\varepsilon>x>0$ and extend it by $0$ to $X$. The magnetic field
$B=dA=\psi'(x)dx\wedge\beta$ has compact support in $X$.

By Theorem \ref{t:thmag}, $\Delta_A$ has purely discrete spectrum with
the Weyl asymptotic law eigenvalues independent of $B$. 

The relative cohomology class $[B]$ lives in the direct sum 
$\oplus_j \partial H^1(M_j,\rz)\subset H^2(\oX,M,\rz)$. The field $B$ is 
non-trapping if at least one of its components in this decomposition 
lives in the image of $H^1(M_j,\zz)$. Since we assume
all $H^1(M_j,\rz)$ to be nonzero, the space of non-trapping magnetic fields 
is a finite union of subspaces of codimension at least $1$.
\qed 

We now show that the cohomological hypothesis about $X$ and $M$
can be satisfied in all dimensions greater than or equal 
to $2$, and different from $3$.

In dimension $2$, take $X=\rz^2$ endowed with the metric
\eqref{pme}. Consider, for instance, the metric $r^{-2p}(dr^2+d\sigma^2)$
given in polar coordinates. Here $M$ is the circle at infinity and
$x=1/r$ for large $r$. Thus $b_1(X)=0$ while $b_1(M)\neq 0$.
The product of this manifold with a closed, connected, simply
connected manifold $Y$ of dimension $k$ yields an example in dimension
$2+k$ with the same properties. Indeed, by the K\"unneth formula,
the first cohomology group
of $\rz^2\times Y$ vanishes, while $H^1(S^1\times Y)\simeq
H^1(S^1)=\zz$. Clearly $k$ cannot be $1$ since the only closed
manifold in dimension $1$ is the circle. Thus the dimension $3$ is
actually exceptional.

For orientable $X$ of dimension $3$, the assumptions $H^1(X)=0$ and
$H^1(M)\neq 0$ cannot be simultaneously fulfilled.
Indeed, we have the following long exact sequence (valid actually
regardless of the dimension of $X$)
\begin{equation} \label{les3}
H^1(\oX)\stackrel{i_M}{\longrightarrow} H^1(M)\stackrel{\delta}{\longrightarrow} H^2(\oX,M).
\end{equation}
If $\dim(X)=3$, the spaces $H^1(\oX)$ and $H^2(\oX,M)$ are isomorphic by
Poincar\'e duality, hence $H^1(\oX)=0$ implies $H^2(\oX,M)=0$ and so
(by exactness) $H^1(M)=0$. It should be possible to build a
non-orientable example in dimension $3$ such that one could apply
Proposition \ref{p:contreex} but we were not able to construct one. 

We finally give an example not covered by Lemma \ref{l:compsupp} but such that
the conclusion of Theorem \ref{t:thmag} holds. We considered so far
magnetic fields on $\overline X$ with vector potential smooth on $\oX$. One
may consider also more singular magnetic fields arising from Proposition
\ref{p:cuspduB}. 
\subsubsection{Example} Let $X$ be any conformally cusp manifold, without
any cohomological assumptions.
Suppose that $B=df\wedge dx/x^2$ where $f$ is a function on $X$ smooth
down to the boundary $M$ of $X$. Assume that $f$ is not constant on any
connected component of
$M$. Then the essential spectrum of the magnetic operator (which is
well-defined by $B$ if $H^1(X)=0$) is empty. This follows from
the fact that $A:=fdx/x^2$ is trapping.

Note that the pull-back to the border of the above magnetic field is zero.

\subsection{The coupling constant effect}\label{s:coupling}
In flat Euclidean space it is shown in \cite{hn}, under some
technical hypotheses, that the spectrum has a limit as the coupling
constant tends to infinity. In contrast, in the next example, we
exhibit the creation of essential spectrum for periodic values of the
coupling constant. We will focus on the properties of $\Delta_{gB}$ 
for some coupling constant $g\in\R$. In order to be able to
exploit the two sides of this work we concentrate here on the metric
\eqref{gp'}. We assume that $M$ is \emph{connected}, $H^1(X,\zz)=0$ and
$H^1(M)\neq 0$.

Given $B$ a magnetic potential with compact support, $g B$ is
non-trapping if and only its class in $H^2(\oX,M;\rz)$ is integral. Let
$G_B$ be the discrete subgroup of those $g\in\rz$ such that
$g B$ is non-trapping. As this subgroup is possibly $\{0\}$,
we start with some exact form $B$ which represents
a nonzero cohomology class in $H^2(\oX,M;\zz)$; then by exactness 
of the relative cohomology long sequence, $[B]$ lives in the image
of the injection $H^1(M;R)\stackrel{\partial}{\longrightarrow}H^2(\oX,M;\rz)$. 
With these restrictions, $G_B$ is a non-zero discrete subgroup of $\qz$.
Now we apply Theorem
\ref{t:thmag} for the trapping case and Theorem \ref{t:mourre0}
for the non-trapping case. We obtain that
\begin{enumerate}
\item For $g\in G_B$, the essential spectrum of
$\Delta_{g B}$ is given by $[\kappa(p), \infty)$, where
$\kappa(p)$ is defined in Proposition \ref{p:thema}. The spectrum of
$\Delta_{g B}$ has no singular continuous part and the
eigenvalues of $\R\setminus \{\kappa(p)\}$ are of finite multiplicity and
can accumulate only in $\kappa(p)$.
\item For $g\notin G_B$, the spectrum $\Delta_{g B}$ is
discrete and if $p\geq 1/n$, the asymptotic of the eigenvalues
depend neither on $g$ nor on $B$.
\end{enumerate} 
We now describe the long-range effect regarding the coupling
constant. Take a state $\phi\in
L^2(X)$ such that $\phi$ is not an eigenvalue of the free Laplacian $\Delta_0$ 
and is located in a energy higher than $\kappa(p)$. Since the Fourier 
transform of an absolutely continuous measure
(comparing to the Lebesgue measure) tends to $0$ at infinity, we
obtain, for each $g\in G_{g B}$ and for each $\cchi$
operator of multiplication by the characteristic function of compact
support that $\cchi e^{it\Delta_{g B}}\phi\rightarrow 0$ as
$t\rightarrow \infty$. If one considers $\cchi$ being 1 above the
support of the magnetic field, then after some time the norm of $\phi$
above this zone is arbitrary small. Classically, the particle stops
interacting with the magnetic field. Let us denote by $\phi'$
the particle at this moment. Now, if we switch on the interaction with
intensity as small as one desires one gets $g\notin G_{g B}$ and then 
the spectrum of $\Delta_{g B}$ is discrete.
Therefore there exists $\cchi'$ operator of multiplication by the
characteristic function of compact support such that $1/T
\int_{0}^{T} \|\cchi' e^{it\Delta_{g B}}\phi'\|^2 dt$ tends to
a positive constant, 
as $T\rightarrow \infty$. The particle is caught
by the magnetic field even thought they are far from being able to
interact classically. 

In other words, switching on the interaction of the magnetic field
with compact support has destroyed the situation of limiting absorption principle.
This is a strong long-range effect.

For the sake of utmost concreteness, take $X=\R^2$ endowed with the metric
$r^{-2p}(dr^2+d\theta^2)$ in polar coordinates, for $r$ big enough
and $1\geq p>0$. The border $M$ is $S^1$ and
$H^2(\overline{X}, M;\R)\simeq \Z$. Then for every closed $2$-form $B$ with compact
support and non-zero integral, the group $G_B$ defined above is non-zero.

\subsection{The case $H^1(X)\neq 0$, the Aharonov-Bohm effect}\label{s:AB} 
Gauge invariance does not hold in this case, so one expects some sort of
Aharonov-Bohm effect \cite{AB}. Indeed, given two vectors potential
arising from the same magnetic field, the associated magnetic Laplacians
$\Delta_A$ and $\Delta_{A'}$ might be unitarily in-equivalent. The
vector potential acquires therefore a certain physical meaning in this case. 

In flat $\R^n$ with holes, some long-range effect appears, for instance
in the scattering matrix like in \cite{R, RY,
RY2, Y}, in an inverse-scattering problem \cite{N,W} or in the
semi-classical regime \cite{BR}. See also \cite{helffer} for the
influence of the obstacle on the bottom of the spectrum.

In all the above cases the essential spectrum remains the same. 
In light of Proposition \ref{p:contreex}, one can expect a much stronger 
effect in our context. We now give some examples of magnetic fields with 
compact support such that there exists a non-trapping vector potential $A$, 
constant in $x$ in a neighborhood of $M$, and a trapping vector
potential $A'$, 
such that $dA=dA'=B$. To ease the presentation, we stick to the metric 
\eqref{gp'}. For $A$, one applies Theorem \ref{t:mourre0} and obtain 
that the essential spectrum of $\Delta_{A}$ is given by $[\kappa(p), \infty)$, 
that $\Delta_{A}$ has no singular continuous part and the eigenvalues 
of $\R\setminus \{\kappa(p)\}$ have finite multiplicity and can accumulate 
only to $\{\kappa(p)\}$. For $A'$, one applies Theorem \ref{t:thmag} to
get the discreteness of the spectrum of $\Delta_{A'}$ and to obtain that 
the asymptotic of eigenvalues depends neither on $A$ nor on $B$ 
for $n\geq 1/p$. In the next section, we describe how generic 
this situation is for hyperbolic manifolds of dimension $2$ and $3$.

The easy step is to construct the non-trapping vector potential $A$ 
constant in a neighborhood of $M$, more precisely we construct $A$ to be 
\emph{maximal} non-trapping, see Definition \ref{def7}. Indeed, prescribe a 
closed $1$-form $\theta$ on $M$ defining integral cohomology $1$-classes on
each component of the boundary, and extend it smoothly to $X$,
constant in $x$ in a neighborhood on $M$, like in the proof of
Proposition \ref{p:contreex}. The magnetic field $B:=dA$ has then
compact support and one may apply Theorem \ref{t:mourre0} for
$\Delta_A$. We now construct $A'$ by adding to $A$ a closed form
$\alpha$, smooth on $\overline{X}$. Since $\alpha$ is closed, 
$A'$ and $A$ define the same magnetic field. 
In light of Remark \ref{addtrap}, $A'$ is trapping if and only if $\alpha$ is. 
We can then apply Theorem \ref{t:thmag} to $\Delta_{A'}$.
It remains to show that closed, trapping $\alpha$ do exist. We start 
with a concrete example.

\begin{example}
Consider the manifold $\oX=(S^1)^{n-1}\times [0,1]$ with a metric $g_p$ as in
\eqref{gp'} near the two boundary components. Let $\theta_i\in \rz$
be variables on the torus $(S^1)^{n-1}$, so $e^{i\theta_i}\in S^1$. Take
the vector potential $A$ to be $0$, it is (maximal) non-trapping. Choose now 
$\alpha=A'$ to be the closed form $\mu d\theta_1$ for some
$\mu\in\rz$. It is constant in a neighborhood of $(S^1)^{n-1}$. The
class $[i_M^*(A')]$ is an integer multiple of $2\pi$ if and only if
$\mu\in\zz$. In other words, $A'$ is non-trapping if and only if 
$\mu\in\rz\setminus\zz$. Note that here the magnetic field $B$ vanishes. 
\end{example}

In order to show the existence of such $\alpha$ in a more general setting, 
we assume that the first Betti number of each connected component of 
the boundary is non-zero. It is enough to find some closed $\alpha$,
smooth on $\oX$, which on each boundary component represents 
a non-zero cohomology class. Then, up to a multiplication by a constant, 
$\alpha$ will be trapping. When $X$ is orientable and $\dim(X)$ 
is $2$ or $3$, one proceeds as follows.

\begin{proposition}\label{p23}
Let $\oX$ be a compact manifold with non-empty boundary $M$.
Assume that one of the following hypotheses holds:
\begin{enumerate}
\item $\dim(X)=2$ and $M$ is disconnected;
\item $\dim(X)=2$ and $X$ is non-orientable;
\item $\dim(X)=3$, $X$ is orientable 
and none of the connected components of $M$ are spheres.
\end{enumerate}
Then there exists a closed smooth form $\alpha\in\Lambda^1(\oX)$, constant in
$x$ near the boundary, such that for all connected components $M_j$ of $M$,
the class $[\alpha_{|M_j}]\in H^1(M_j;\rz)$ is non-zero.
\end{proposition}
\proof
Any cohomology class on $\oX$ admits a smooth representative $\alpha$ 
up to the boundary. Moreover since $\alpha$ is closed, one can choose
$\alpha$ to be constant in $x$ near the boundary using Lemma
\ref{l:compsupp}. Thus, in cohomological terms, the proposition is
equivalent to finding a class $[\alpha]\in H^1(\oX)$ whose 
pull-back to each connected component of $M$ is non-zero, i.e.,
$H^1(M_j)\ni i_{M_j}[\alpha]\neq 0$. 

Consider first the case $\dim(X)=2$. If $X$ is 
non-orientable, $H^2(\oX,M)=0$ so $\delta$ is the zero map. 
If $X$ is oriented, the boundary 
components are all oriented circles, so $H^1(M_j)\simeq \rz$; the compactly
supported cohomology $H^2(\oX,M)$ is isomorphic to $\rz$ via the integration map,
and the boundary map $\delta:H^1(M)\to H^2(\oX,M)$ restricted to $H^1(M_j)$
is just the identity map of $\rz$ under these identifications. Thus the kernel
of $\delta$ is made of $\nu$-tuples (where $\nu$ is the number of boundary
components) $(a_1,\ldots,a_\nu)$ of real numbers, with the
constraint $\sum a_j=0$ in the orientable case. By exactness,
this space is also the image of the restriction map $H^1(\oX)\to H^1(M)$.
Clearly there exist such tuples with non-zero entries, provided
$\nu\geq 2$ in the orientable case. Thus the conclusion follows for $\dim(X)=2$.

Assume now that $\dim(X)=3$. Then the maps $i_M$ and $\delta$ from the relative long
exact sequence \eqref{les3} are dual to each other under the intersection
pairing on $M$, respectively on $H^1(\oX)\times H^2(\oX,M)$:
\[\int_M i_M(\alpha)\wedge\beta=\int_X \alpha\wedge\delta \beta.\]
These bilinear pairings are non-degenerate by Poincar\'e duality; in particular
since the pairing on $H^1(M)$ is skew-symmetric, it defines a symplectic form.
It follows easily that the subspace
\[L:=i_M(H^1(X))\subset H^1(M)\]
is a Lagrangian subspace (i.e., it is a maximal isotropic subspace for the
symplectic form). Now the symplectic vector space $H^1(M)$ splits into the
direct sum of symplectic vector spaces $H^1(M_j)$, By hypothesis, 
the genus $g_j$ of the oriented surface $M_j$ is at least $1$ so $H^1(M_j)$ 
is non-zero for all $j$. It is clear that the projection of $L$ on 
each $H^1(M_j)$ 
must be non-zero, otherwise $L$ would not be maximal. Hence, there exists an
element of $L=i_M(H^1(X))$ which restrict to non-zero classes in each
$H^1(M_j)$, as desired.
\qed

\begin{remark}\label{r:int}
If one is interested in some coupling constant effect, it is interesting 
to choose the closed form $\alpha$ so that every 
$[\alpha_{|M_j}]$ are non-zero integral classes. In dimension $2$, 
this amounts to choosing non-zero integers 
with zero sum. In dimension $3$, as $L=i_M(H^1(X))=i_M(H^1(X,\zz))\otimes \rz$
is spanned by integer classes, we can find an 
integer class in $L$ with non-zero projection on all $H^1(M_j)$. For real $g$, 
the vector potential $g\alpha$ is therefore non-trapping precisely for $g$ 
in a discrete subgroup $g_0\zz$ for some $g_0\in\qz$.
\end{remark}

\section{Application to hyperbolic manifolds}\label{hyp}

We now examine in more detail how this Aharonov-Bohm effect arises in the context 
of hyperbolic manifolds of finite volume in dimension $2$ and $3$. 

These are conformally cusp manifold with $p=1$, with unperturbed
metric of the form \eqref{mc} and such that every component $M_j$ of the
boundary is a circle when $\dim(X)=2$, respectively a flat torus when $\dim(X)=3$:
indeed, outside a compact set, the metric takes the form
$g=dt^2+e^{-2t}h$, where $t\in[0,\infty)$, and this is of the form \eqref{mc} after the change
of variables $x:=e^{-t}$.
We denote by $\oX$ the compactification of $X$ by requiring that $x:=e^{-t}$ be a
boundary-defining function for ``infinity''.
These manifolds and their boundary components always have non-zero first
Betti number.

For complete hyperbolic \emph{surfaces} with cusps, every smooth $2$-form $B$ 
on $\oX$ must be exact because $H^2(\oX;\rz)$ is 
always zero for a non-closed surface. Call $A$ a smooth primitive of $B$. 
If $B$ vanishes at $M$ then $A$ is necessarily closed over $M$. 
In terms of cohomology classes, we have 
$[B]=\delta_\rz [A_{|M}]$ where $\delta_\rz$ is the connecting morphism from the sequence
\eqref{les3} with real coefficients. Notice that $A$ can be 
chosen to define an integer class on $M$ if and only if $[B]$ is integer. 
Indeed, $H^2(\oX;\zz)$ is also $0$ for a compact surface with non-empty 
boundary, so if $[B]$ is integer, it must lie in the image of $\delta_\zz$
for the sequence \eqref{les3} with integer coefficients. Conversely, if $[A]$ is 
integer, i.e., $[A]=2\pi[A_\zz]$ (see diagram \eqref{cohoz} with $M$ in the place of 
$M_\cB$, where the horizontal maps are now \emph{surjective}) then
$[B]=2\pi\delta_\zz [A_\zz]$ is also integer. We summarize these
remarks in the following 
\begin{corollary}\label{81}
Let $B$ be a smooth $2$-form on the compactification of a complete
hyperbolic surface $X$ with cusps, and denote by $[B]\in H^2(\oX,M,\rz)$
its relative cohomology class.
\begin{itemize}
\item If either $X$ is non-orientable, or $X$ is orientable with at least 
two cusps, then $B$ admits both trapping and non-trapping vector potentials.
\item If $X$ is orientable with precisely one cusp, then $B$ admits only 
non-trapping vector potentials if $[B]$ is integral, while if $[B]$ is not integral
then $B$ admits only trapping vector potentials.
\end{itemize} 
\end{corollary}
\begin{proof}
First note that $B$ is closed since it is of maximal degree; it is 
exact since $H^2(\oX)=0$ for every surface with boundary; moreover 
its pull-back to the $1$-dimensional boundary also vanishes, 
so $B$ defines a relative de Rham class. By Corollary \ref{cohcl}, 
the existence of trapping and non-trapping vector potentials 
depends only on this class.

If $X$ is orientable and has precisely $1$ cusp, then the map 
$\delta:H^1(M)\to H^2(\oX,M)$ is an isomorphism both for real and for 
integer coefficients. Thus $[B]$ is integer if and only if 
$[A_{|M}]$ is integer. Since the boundary is connected, $A$ is trapping 
if and only if the cohomology class of its restriction to the boundary 
is non-integer.

If $X$ is oriented and has at least two cusps, identify $H^2(\oX,M)$ and each
$H^1(M_j)$ with $\zz$, so that the boundary map restricted to
$H^1(M_j)$ is the identity. We can write $[B]$ first as a sum $\sum\alpha_j$ of 
non-integer numbers, then also as a sum where at least one term is integer. Let
$A$ be a $1$-form on $\oX$ which restricts to closed forms of cohomology class
$\alpha_j$ on $M_j=S^1$. Then $B-dA$ represents the $0$ class in $H^2(\oX,M)$, so 
after adding to $A$ a form vanishing at the boundary, we can assume that
$B=dA$. Now when all $\alpha_j$ are non-integers, $A$ is trapping, while 
in the other case it is non-trapping as claimed.

If $X$ is non-orientable, the class $[B]$ vanishes. It is enough to find trapping 
and non-trapping vector potentials for the zero magnetic field, which is done 
as in the orientable case.
\end{proof}

When $\dim(X)=3$, we have:

\begin{corollary}
Let $X$ be an orientable complete hyperbolic $3$-manifold of finite volume. 
Then every magnetic field $B$ smooth on the compactification 
$\oX$ admits trapping vector potentials.

Assume that the pull-back of $B$ to the boundary $M$ vanishes.
If $X$ has precisely one cusp, then there exists a rational 
(i.e., containing integer classes) 
infinite cyclic subgroup $G\subset H^2(\oX,M,\rz)$ 
so that $B$ admits a non-trapping vector potential if and only if $[B]\in G$.
In general, one of the following alternative statements holds:
\begin{enumerate}
\item Either every magnetic field smooth on $\oX$ and vanishing at $M$
admits a non-trapping vector potential, or
\item Generically, magnetic fields smooth on $\oX$ and vanishing at $M$ do 
not admit non-trapping vector potentials.
\end{enumerate}
There exists moreover $q\in \zz^*$ such that if 
$[B]$ is integer, then $qB$ admits non-trapping vector potentials.
\end{corollary}

\begin{proof}
For the existence of trapping vector potentials we use the closed form
$\alpha$ from Proposition \ref{p23}. Let $A$ be any vector potential for $B$. 
It suffices to note that for $u\in\rz$,
the form $A+u\alpha$ is another vector potential for $B$, which is trapping
on each connected component of $M$ except possibly for 
some discrete values of $u$.

Let $h$ denote the number of cusps of $X$. Both the Lagrangian subspace
$L\subset H^1(M)$ and the image space $\partial H^1(M)\subset H^2(\oX,M)$ 
have dimension $h$. By hypothesis, the cohomology class of $B$ on $\oX$ is $0$
so by exactness of \eqref{cohoz},
the relative cohomology class $[B]$ lives in $\partial H^1(M)$.

Assume first that $X$ has precisely one cusp.
Let $A$ be a vector potential for $B$ (smooth on $\oX)$. 
We can change $A$ by adding to it any class in the line $L$ 
without changing $[B]$.
Notice that $H^1(M)=\zz^2$. The line $L$ has an integer generator 
(given by the image of $H^1(\oX,\zz)\to H^1(M,\zz)$). Without loss of
generality, we can assume that $L$ is not the horizontal axis in $\zz^2$. It
follows that the translates of all integer points in $\zz^2$ in
directions parallel to $L$ form a discrete subgroup of $\qz$. Thus $B$
admits non-trapping vector potentials if and only if the cohomology
class $[B]$ inside the $1$-dimensional image $\partial H^1(M)$ lives
inside a certain infinite cyclic 
discrete subgroup. In particular, if $B$ is irrational
(i.e., no positive integer multiple of $B$ is an integral class) then
$B$ does not admit non-trapping vector potentials.

In the general case, assume first that there exists a boundary component $M_j$ 
so that $L$ projects surjectively onto $H^1(M_j,\rz)$. Let $A$ be a vector 
potential for an arbitrary magnetic field $B$ which vanishes at $M$.
Let $[A'_j]\in H^1(M_j,\rz) $ be such that 
$[A]_{H^1(M_j)}+[A'_j]$ is integer. Let $[A']\in L$ be an element whose 
component in $H^1(M_j)$ is $[A'_j]$. Choose a representative $A'$ and
extend it to a smooth $1$-form on $\oX$, constant in $x$ near the boundary.
Then $A+A'$ is a non-trapping vector potential. From the definition of 
$L$, the form $dA'$ defines the zero class in 
relative cohomology, so from Corollary \ref{cohcl} we get the assertion on $B$.

If the assumption on $L$ is not fulfilled, we claim that 
\[\hat{L}_j:=L\cap \oplus_{i\neq j}H^1(M_i)\] 
has dimension $h-1$ for all $j$.
Indeed, this dimension cannot be $h$ (since the projection of $L$
on $H^1(M_j)$ is not zero) and it cannot be $h-2$ (since the projection
is not surjective). It follows easily that $\hat{L}_j$ is a Lagrangian subspace
of $\oplus_{i\neq j}H^1(M_i)$. Let $v$ be a vector in $L\setminus \hat{L}_j$.
The component $\hat{v}_j$ of $v$ in $\oplus_{i\neq j}H^1(M_i)$ is 
clearly orthogonal (with respect to the symplectic form) to $\hat{L}_j$, so
by maximality it must belong to $\hat{L}_j$. Thus we may subtract 
this component to obtain, for each $j$, a non-zero element of $L\cap
H^1(M_j)$. 
These elements may be taken integral since $L$ has integer generators. Since
$L=\ker\partial$, it follows that the image of $\partial$ is the direct sum 
of the images $\partial(H^1(M_j))$.
As in the case of only one cusp, we see that $B$ has a non-trapping potential
if and only if at least one of the components of $[B]$ 
in this decomposition belong to a certain cyclic subgroup containing integer 
classes.

Set $q$ to be the least common denominator of the generators of these 
subgroups for all $j$. If $[B]$ is integer, it follows that every 
vector potential for $qB$ must be maximal non-trapping.
\end{proof}

\appendix

\section{The $C^1$ condition in the Mourre theory}\label{s:c1} 
\renewcommand{\theequation}{A.\arabic{equation}}
\setcounter{equation}{0}
In this appendix, we give a general criterion of its own interest to 
check the, 
somehow abstract, hypothesis of regularity $\cC^1$ which is a key
notion in the Virial theorem within Mourre's theory, see \cite{ABG}
and \cite{GG0}. Let $A$ and $H$ be two self-adjoint operators in a
Hilbert space $\Hr$. The commutator $[H,iA]$ is defined in the sense of forms
on $\Dc(A)\cap\Dc(H)$. Suppose that the commutator $[H,iA]$ extends to
$\Bc(\Dc(H), \Dc(H)^*)$ and denote by $[H,iA]_0$ the extension. Suppose also
that the following Mourre estimate holds true on an open interval $\cI$, 
i.e.\ there is a constant $c>0$ and a compact operator $K$ such that
\begin{equation}\label{e:mourreK}
E_\cI(H)[H,iA]_0 E_\cI(H)\geq c E_\cI(H) +K,
\end{equation} 
where $E_\cI(H)$ denotes the spectral measure of $H$ above $\cI$. 
Take now $\lambda\in\cI$ which is not an eigenvalue of $H$. Set 
$\cI_n:=(\lambda-1/n, \lambda+1/n)$.
Then $E_{\cI_n}(H)$ tends strongly to $0$ as $n\to\infty$, so 
$E_{\cI_n}(H)KE_{\cI_n}(H)$ 
tends in norm to 
$0$. Hence for $n$ big enough and for some $0<c'\leq c$, 
one gets the strict Mourre estimate
\begin{equation}\label{e:mourrestrict}
E_{\cI_n}(H)[H,iA]_0 E_{\cI_n}(H)\geq c' E_{\cI_n}(H).
\end{equation} 

By supposing that $H\in\cC^1(A)$ (see below) or that 
$e^{itA}\Dc(H)\subset\Dc(H)$, the Virial theorem holds true, i.e.\ 
$\langle f, [H,iA]_0 f\rangle=0$ for every eigenvector $f$ of $H$. Note that 
$f$ has no reason to lie in $\Dc(A)$ and that the expansion of the 
commutator $[H-\lambda, iA]$ over $f$ is formal. 

The Virial theorem is crucial to the study of embedded eigenvalues of $H$. 
Assuming \eqref{e:mourreK}, it implies the local finiteness of the 
point spectrum of $H$ over $\cI$, i.e., that the sum of the multiplicities 
of the eigenvalues of $H$ inside $\cI$ is finite. 
To see this, apply \eqref{e:mourreK} to a infinite 
sequence $(f_n)_{n\in\N}$ of orthonormal eigenvectors of $H$. 
Then, since $\langle f_n, K f_n \rangle$ tends to $0$ as $n$ goes to 
infinity, one obtains a contradiction with the positivity of $c$. 
Assuming \eqref{e:mourrestrict}, the Virial theorem implies directly
that $H$ has no eigenvalue in $\cI_n$.

We stress that the hypothesis $[H,iA]_0\in\Bc(\Dc(H), \Dc(H)^*)$ 
does not imply the Virial theorem. A counterexample is given in 
\cite{GG0}. If one adds some conditions on the second order commutator of 
$H$ and $A$, e.g.\ like in \cite{CFKS}, one deduces from \eqref{e:mourrestrict}
a limiting absorption principle and therefore the absence of eigenvalues in 
$\cI_n$. In turn, assuming \eqref{e:mourreK}, we deduce that the
set of eigenvalues of $H$ in $\cI$ is closed. 
It is not known whether the multiplicity of the 
point spectrum must be locally finite when the Virial theorem does not hold. 
 
Checking the Virial theorem, or a sufficient condition for it 
like the $C^1$ condition, is sometimes omitted in the Mourre analysis 
in a manifold context. To our knowledge, no result exists actually to
show directly the $C^1$ regularity in a manifold context. On a class
of exponentially growing manifolds, Bouclet \cite{bouclet} 
circumvents the problem by showing a stronger fact, i.e., 
the invariance of the domain. This method does not seem to work for our 
local conjugate operator, see Section \ref{s:conjop}. Besides
giving an abstract criterion for the $\cC^1$ condition, 
we will explain under which additional condition we can 
recover the invariance of the domain from it.
 
Given $z\in\rho(H)$, we denote by $R(z)=(H-z)^{-1}$. For $k\in\N$, we recall
that $H\in\cC^k(A)$ if for one $z\notin\sigma(H)$ (then for all
$z\notin\sigma(H)$) the map $t\mapsto e^{-itA}R(z)e^{itA}$ is $C^k$
in the strong topology. 
We recall a result following from Lemma 6.2.9 and Theorem 6.2.10 of 
\cite{ABG}. 
\begin{theorem}\label{th:abg} 
Let $A$ and $H$ be two self-adjoint operators in the Hilbert space
$\Hr$. The following points are equivalent: 
\begin{enumerate} 
\item $H\in\cC^1(A)$. 
\item For one (then for all) $z\notin \sigma(H)$, there is a finite
$c$ such that 
\begin{align}\label{e:C1a} 
|\langle A f, R(z) f\rangle - \langle R(\overline{z}) f, Af\rangle| \leq c 
\|f\|^2, \mbox{ for all $f\in\Dc(A)$}. 
\end{align} 
\item 
\begin{enumerate} 
\item [a.]There is a finite $c$ such that for all $f\in \Dc(A)\cap\Dc(H)$: 
\begin{equation}\label{e:C1b} 
|\langle Af, H f\rangle- \langle H f, Af\rangle|\leq 
 c(\|H f\|^2+\|f\|^2). 
\end{equation} 
\item [b.] For some (then for all) $z\notin \sigma(H)$, the set
$\{f\in\Dc(A) \mid R(z)f\in\Dc(A)$ and $R(\overline{z})f\in\Dc(A)
\}$ is a core for $A$. 
\end{enumerate} 
\end{enumerate} 
\end{theorem} 
Note that in practice, condition (3.a) is usually easy to check
and follows from
the construction of the conjugate operator. The condition (3.b) could be 
more delicate. This is addressed in the next lemma, 
inspired by \cite{bouclet}.
 
\begin{lemma}\label{l:C1} 
Let $\Dr$ be a subspace of $\Hr$ such that $\Dr\subset\Dc(H)\cap\Dc(A)$, 
$\Dr$ is a core for $A$ and $H\Dr\subset\Dr$. Let $(\cchi_n)_{n\in
\N}$ be a family of bounded operators such that
\begin{enumerate} 
\item $\cchi_n\Dr\subset \Dr$, $\cchi_n$ tends strongly to $1$
as $n\to\infty$,
and $\sup_n \|\cchi_n\|_{\Dc(H)}<\infty$.
\item $A\cchi_n f \rightarrow Af$,
for all $f\in\Dr$, as $n\to\infty$, 
\item There is $z\notin \sigma(H)$, such that $\cchi_n R(z) \Dr \subset 
 \Dr \mbox{ and } \cchi_n R(\overline{z}) \Dr \subset\Dr$.
\end{enumerate} 
Suppose also that for all $f\in\Dr$
\begin{equation}\label{e:comm} 
\lim_{n\rightarrow \infty} A [H, \cchi_n] R(z)f= 0
\mbox{ and } 
\lim_{n\rightarrow \infty} A [H, \cchi_n] R(\overline{z})f=
0.
\end{equation} 
Finally, suppose that there is a finite $c$ such that 
\begin{equation}\label{e:C1} 
\hspace*{1cm} |\langle Af, H f\rangle- \langle H f, Af\rangle|\leq 
 c(\|H f\|^2+\|f\|^2),\quad \forall f\in\Dr. 
\end{equation} 
Then one has $H\in\cC^1(A)$. 
\end{lemma} 
Note that \eqref{e:comm} is well defined by expanding the commutator 
$[H,\cchi_n]$ and by using (3) and $H\Dr\subset \Dr$. 
\proof 
By polarization and by applying \eqref{e:C1} to $\cchi_n 
R(\overline{z})f$ and to $\cchi_n R(z)f$, with $f\in\Dr$,
we see that there exists $c<\infty$ such that 
\begin{align} \nonumber
\lefteqn{ |\langle A\cchi_n R(\overline{z})f, H \cchi_n R({z})f\rangle- 
\langle H \cchi_n R(\overline{z}) f, A\cchi_n R({z}) f\rangle|}\\&&
\leq c\|(H+i)\cchi_n R(\overline{z})f\|\cdot \|(H+i)\cchi_n R(z)f\|, \label{e:A7}
\end{align} 
for all $f\in\Dr$. By condition (1), the right-hand side is bounded by
$C\|f\|^2$ for some $C$. 
We expand the left hand side of \eqref{e:A7} by commuting $H$ with $\cchi_n$: 
\begin{align*} 
& | \langle \cchi_n R(\overline{z})f, A\cchi_n f\rangle- \langle 
 A\cchi_n f, \cchi_n R({z}) f\rangle\\ 
& + \langle \cchi_n R(\overline{z})f, A[H,\cchi_n] R({z})f\rangle- 
\langle A[H, \cchi_n] R(\overline{z}) f, \cchi_n R({z}) f\rangle|. 
\end{align*} 
Using \eqref{e:comm}, the second line vanishes as $n$ goes to infinity. 
Taking in account the assumptions (1) and (2), we deduce:
\begin{equation*} 
|\langle R(\overline{z})f, A f\rangle- \langle A f, R({z}) f\rangle | 
\leq C\|f\|^2, \quad \forall f\in\Dr. 
\end{equation*} 
Finally, since $\Dr$ is a core for $A$, we obtain 
\eqref{e:C1a}. We conclude that $H\in\cC^1(A)$. \qed 
 
The hypotheses of the lemma are easily satisfied in a manifold context
with $H$ being the Laplacian and $A$ its
conjugate operator, constructed as a localization on the ends of the
generator of dilatation like for instance in \cite{FH}. Let $\Dr=\cunc(X)$ 
and $\cchi_n$ a family of operators of multiplication by smooth
cut-off functions with compact support. The fact that $A$ is self-adjoint
comes usually by some consideration of $C_0$-group associated to some vector
fields and using the Nelson Lemma and the invariance of $\Dr$ under the
$C_0$-group give that $\Dr$ is a core for $A$, see Remark \ref{r:sa}. The
hypothesis (3) follows then by elliptic regularity. The
only point to really check is \eqref{e:comm}. At this point one needs to
choose more carefully the family $\cchi_n$. In Lemma
\ref{l:regu0}, we show that the hypotheses of Lemma \ref{l:C1}
hold for the standard conjugate operator and for our local conjugate operator. 

The invariance of the domain is desirable in order to deal, 
in a more convenient way, directly with operators and no longer
with resolvents. On a manifold the $C_0$-group $e^{itA}$ is not explicit
and it could be delicate to deal with the domain of $H$ directly. 
However, one may obtain this invariance of the domain using
\cite{GG0} and a $C^1(A)$ condition. We recall: 
 
\begin{lemma} If $H\in \cC^1(A)$ and $[H, iA]:\Dc(H)\rightarrow
\Hr$ then $e^{itA}\Dc(H)\subset\Dc(H)$, for all
$t\in\R$. 
\end{lemma} 

In light of this lemma, one understands better the importance of having
some $C^1(A)$ criteria. On one hand, one can easily check the
invariance of the domain. On the other hand, if the commutator
belongs only to $\Bc(\Dc(H), \Dc(H)^*)$ and not $\Bc(\Dc(H), \Hr)$,
one may turn to another version of the Mourre Theory like in
\cite{ABG} when $H$ has a spectral gap or like in \cite{GJ, S} in the
other case.

\section{Finite multiplicity of $L^2$ eigenvalues}\label{fm}

In order to classify all maximal symmetric extension of a given
cusp-elliptic operator $H$ (see section \ref{s:gene} for definitions), one 
computes the defect indices, i.e., $\dim\ker (H^*\pm i)$. If they are equal 
and finite, one concludes that
all maximal symmetric extensions of $H$ are self-adjoint. By the Krein
formulae, one hence obtains that the difference of the resolvent of two 
maximal
extensions is finite rank. This implies by Weyl's theorem that the essential
spectrum is the same for all self-adjoint extensions. Moreover, by Birman's
theorem, the wave operators associated to a pair of such extensions exist
and are complete. On the other hand, note that if the defect indices 
are not finite, one may have maximal symmetric extensions which are not 
self-adjoint. 

It is also interesting to control the multiplicity of eigenvalues embedded in
the essential spectrum. 

In the next lemma we assume that $X$ is a conformally cusp manifold 
with respect to the metric \eqref{pme}. We fix a vector bundle $E$ over $\oX$ 
(for instance the bundle of cups differential forms, although in this paper 
we only use the case where $E$ is the trivial bundle $\cz$) endowed with a smooth 
metric up to $\partial \oX=M$.
 
\begin{lemma}\label{l:indices}
Let $\Delta$, acting on $\cun_c(X, E)$, be a cusp-elliptic differential 
operator in $x^{-2p}\Diff^k(X,E)$ for some $p,k>0$. Then the dimension of any
$L^2$-eigenspace of $\Delta^*$ is finite.
\end{lemma}

Remember that if the operator $\Delta$ is bounded from below, then the defect
indices are the same. This lemma guarantees that they are also finite. This
point is not obvious when the manifold is not complete even if $\Delta$
is a Laplacian. Of course,
this result is based on ideas that can be traced back to \cite{melaps}
and which are today quite standard. This lemma generalizes a result of
\cite{GMo}. 

\proof 
We start by noticing that $\Delta$ can be regarded as an unbounded
operator in a larger $L^2$ space. Namely, let $L^2_{\varepsilon}$ be
the completion of $\cunc(X,E)$ with respect to the volume form
$e^{-\frac{2\varepsilon}{x}}dg_p$ for some $\varepsilon>0$. Clearly then
$L^2_{\varepsilon}$ contains $L^2$. A distributional solution of
$\Delta-\lambda$ in $L^2$ is evidently also a distributional
solution of $\Delta-\lambda$ in $L^2_{\varepsilon}$. Thus the conclusion
will follow by showing that $\Delta$ has in $L^2_{\varepsilon}$
a unique closed extension with purely discrete spectrum. The strategy
for this is by now clear. First we conjugate $\Delta$ through the isometry
\begin{align*}
L^2_{\varepsilon}\to L^2,&&\phi\mapsto e^{-\frac{\varepsilon}{x}}\phi.
\end{align*}
We get an unbounded operator $e^{-\frac{\varepsilon}{x}} 
\Delta e^{\frac{\varepsilon}{x}}$
in $L^2$, which is unitarily equivalent to $\Delta$ 
(acting in $L^2_{\varepsilon}$).
Essentially from the definition, see Section \ref{s:gene},
\[\cN(x^{2p}e^{-\frac{\varepsilon}{x}}\Delta 
e^{\frac{\varepsilon}{x}})(\xi)= \cN(x^{2p}\Delta)(\xi+i\varepsilon).\] 
The normal operator is a polynomial in $\xi$, in particular it is entire. 
Then by analytic Fredholm theory \cite[Prop.\ 5.3]{melaps}, the set of 
complex values of $\xi$ for which $\cN(x^{2p}\Delta)(\xi)$ 
is not invertible, is discrete. Thus there exists $\varepsilon \pg 0$ 
such that $\cN(x^{2p}\Delta)(\xi+i\varepsilon)$ is invertible for 
all $\xi\in\rz$. For such $\varepsilon$ the operator $\Delta$ in 
$L^2_\varepsilon$ is unitarily equivalent to a fully elliptic cusp operator 
of order $(k,2p)$ in $L^2$. It is then a general fact about
the cusp algebra \cite[Theorem 17]{wlom} that such an operator 
has a unique closed extension and admits a compact inverse modulo 
compact operators in $L^2_\varepsilon$. 
In particular, its eigenvalues have finite multiplicity. 
As noted above, the eigenspaces of $\Delta$ in $L^2$ are contained in the 
eigenspaces of $\Delta$ in $L^2_\varepsilon$ for the same eigenvalue.
\qed

As a corollary, the magnetic Laplacians for the metric \eqref{pme}
and for vector potentials \eqref{aft} which are smooth cusp $1$-forms, 
have finite multiplicity eigenvalues.
 
\section{Stability of the essential spectrum} \label{s:stab}
It is well-known that the essential spectrum of an elliptic differential
operator on a complete manifold can be computed by cutting out a compact
part and studying the Dirichlet extension of the remaining operator on the 
non-compact part (see e.g., \cite{donelili}). This result is obvious using 
Zhislin sequences, but the approach from loc.\ cit.\ fails 
in the non-complete case. For completeness, we give below a proof
which has the advantage to hold in a wider context and for a wider
class of operator, pseudodifferential operators for instance. 

We start with a general lemma.
We recall that a \emph{Weyl sequence} for a couple $(H,\lambda)$ with $H$ a
self-adjoint operator and $\lambda\in\rz$, is a sequence
$\varphi_n\in\Dc(H)$ such that $\|\varphi_n\|=1$,
$\varphi_n\rightharpoonup 0$ (weakly) and such that
$(H-\lambda)\varphi_n\rightarrow 0$, as $n$ goes to infinity.
It is well-known that $\lambda\in\sigma_{\rm ess}(H)$ if and only if
there is Weyl sequence for $(H,\lambda)$.

\begin{lemma}\label{l:cut}
Let $H$ be a self-adjoint operator in a Hilbert space $\Hr$.
Let $\varphi_n$ be a Weyl sequence for the couple $(H,\lambda)$.
Suppose that there is a closed operator $\Phi$ in $\Hr$ such that:
\begin{enumerate}
\item $\Phi \Dc(H)\subset\Dc(H)$,\item $\Phi(H+i)^{-1}$ is compact,
\item $[H,\Phi]$ is a compact operator from $\Dc(H)$ to $\Hr$. 
\end{enumerate}
Then there is $\widetilde{\varphi}_n\in\Dc(H)$ such that
$(1-\Phi)\widetilde{\varphi}_n$ is a Weyl sequence for $(H,\lambda)$.
\end{lemma}
\proof
First we note (2) implies that $\Phi\varphi_n$ goes to $0$. Indeed, we have
$\Phi\varphi_n=\Phi(H+i)^{-1}\big((H-\lambda)\varphi_n 
+(i+\lambda) \varphi_n \big)$ and the bracket goes weakly to $0$. Similarly, 
using (3) we get that $[H,\Phi]\varphi_n\rightarrow 0$. Therefore we obtain 
$\|(1-\Phi)\varphi_n\|\geq 1/2$ for $n$ large enough. We set 
$\widetilde\varphi_n:=\varphi_n/\|(1-\Phi)\varphi_n\|$. Note that 
$(1-\Phi)\widetilde\varphi_n\rightharpoonup 0$. Finally, 
$(H-\lambda)(1-\Phi)\widetilde\varphi_n\rightarrow 0$ since 
$[H,\Phi]\varphi_n\rightarrow 0$.
\qed

This shows that the essential spectrum is given by a ``non-compact''
part of the space. We now focus on Friedrichs extension. Given a dense
subspace $\Dr$ of a Hilbert space $\Hr$ and a positive symmetric
operator on $\Dr$. Let $\Hr_1$ be the completion of $\Dr$ under the
norm given by $\Qr(\varphi)^2=\langle H\varphi,
\varphi\rangle+\|\varphi\|^2$. 
The domain of the Friedrichs extension of $H$, is
given by $\Dc(H_\Fr)=\{f\in \Hr_1 \mid \Dr\ni g\mapsto \langle
Hg,f\rangle+\langle g,f\rangle $ extends to a norm continuous function
on $\Hr\}$. For each $f\in\Dc(H_\Fr)$, there is a unique $u_f$ such
that $\langle 
Hg,f\rangle + \langle g,f\rangle = \langle g,u_f\rangle$, by Riesz
theorem. 
The Friedrichs extension of $H$ is defined by setting $H_\Fr
f:=u_f-f$. It is a self-adjoint extension of $H$, see \cite{RS2}.

Let $(X,g)$ be a smooth Riemannian with distance $d$. We fix $K$ a
smooth compact sub-manifold of $X$ of same dimension. We endow it with
the induced Riemannian metric. We set $X'=X\setminus K$. In the
following, we embed $L^2(X')$ in $L^2(X)$. We will need the next 
definition within the proof. 
\begin{definition}
 We say that $\Phi$ is a (smooth) \emph{cut-off
function} for $K$ if $\Phi\in\cC_c^\infty(X)$ and $\Phi|_K=1$. 
We say that it is an $\varepsilon$-cut-off is $\supp(\Phi)\subset
B(K,\varepsilon)$.
\end{definition} 
We are now able to give a result of stability of the essential
spectrum. 

\begin{proposition}\label{p:cut}
Let $d$ be a differential form of order $1$ on
$C^\infty_c(X)\rightarrow C^\infty(X,\Lambda^1)$ with injective
symbol away from the $0$ section of the cotangent bundle. We denote by
$d_X$ and $d_{X'}$ the closure of $d$ in $L^2(X)$ and $L^2(X')$,
respectively. Consider $\Delta_X:=d_X^*d_X$ and
$\Delta_{X'}=d_{X'}^*d_{X'}$, the Friedrichs extensions of the
operator $d^*d$, acting on $\cC_c^\infty(X)$ and $\cC_c^\infty(X')$,
respectively. One has $\sigma_{\rm{ess}}(\Delta_X)=
\sigma_{\rm{ess}}(\Delta_{X'})$. 
\end{proposition} 

\proof
Let $f\in L^2(B(K, \varepsilon)^c)$ and let $\Phi$ be a
$\varepsilon$-cut-off for $K$. We first show that $f\in\Dc(d_X)$ 
if and only if $f\in\Dc(d_{X'})$. Suppose that $f\in\Dc(d_X)$, then
for all $\eta>0$, there is $\varphi\in \cC_c^\infty(X)$ such that
$\|f-\varphi_n\|+ \|df-d\varphi_n\|<\eta$. Because of the support of
$f$, one obtain that $\|\Phi\varphi\|<\eta\|\Phi\|_\infty$ and that
$\|[d,\Phi]\varphi\|<\eta\|[d,\Phi]\|_\infty$. Therefore
$(1-\Phi)\varphi_n\in \cC_c^\infty(X')$ and is Cauchy in 
$\Dc(d_{X'})$, endow with the graph norm. By uniqueness of the limit, one
obtains that $f\in\Dc(d_{X'})$ and that $d_X f=d_{X'}f$. The opposite
implication is obvious. Using again the $\varepsilon$-cut-off, one shows that
$g\in\Dc(d^*_X)$ if and only if $g\in\Dc(d^*_{X'})$ and that
$d^*_{X}g= d^*_{X'}g$ for $g\in L^2(\Lambda^1 (B(K, \varepsilon)^c)$.
Finally, we obtain that $f\in\Dc(\Delta_X)$ if and only if
$f\in\Dc(\Delta_{X'})$ and that $\Delta_X f=\Delta_{X'} f$, for $f\in L^2(B(K,
\varepsilon)^c)$. 
 
From the definition of the Friedrichs extension and the injectivity of
the symbol of $d$, the domain of $\Delta_X$, $\Delta_{X'}$
is contained in $H^1_0(X)\cap H^2_\loc(X)$, respectively in
$H^1_0(X')\cap H^2_\loc(X')$. By taking the same $\Phi$ as above and
using the Rellich-Kondrakov lemma, the hypotheses of Lemma \ref{l:cut}
are satisfied. Finally, we apply it to $\Delta_X$ and $\Delta_{X'}$
and since the Weyl sequence is with support away from $K$, the
first part of the proof gives us the double inclusion of the essential
spectra. \qed 
 
\bibliographystyle{plain}

Sylvain Gol\'enia, 

Mathematisches Institut der Universit\"at Erlangen-N\"urnberg
Bismarckstr.\ 1 1/2, 91054 Erlangen, Germany

E-mail: golenia@mi.uni-erlangen.de

Sergiu Moroianu, 

Institutul de Matematic\u{a} al Academiei Rom\^{a}ne,
P.O. Box 1-764, RO-014700 Bucha\-rest, Romania 

and

\c Scoala Normal\u a Superioar\u a Bucharest, 
Calea Grivi\c tei 21, Bucharest, Romania

E-mail: moroianu@alum.mit.edu 

\end{document}